\documentclass[12pt,a4paper]{article}

\usepackage{amsmath}
\usepackage{amsfonts}
\usepackage{amssymb}

\newtheorem{defi}{Definition}
\newtheorem{thm}{Theorem}
\newtheorem{lem}{Lemma}
\newtheorem{prop}{Proposition}

\newtheorem{defprop}{Definition \& Proposition}

\newcommand*{\C}{\mathbb{C}}%.............................C
\newcommand*{\Q}{\mathbb{Q}}%.............................Q
\begin{document}
\title{\bf Zeta Functions for $G_2$ and Their Zeros~}
\author{Masatoshi SUZUKI\ \ \ {\footnotesize{and}}\ \ \  Lin WENG\,\thanks{MS is fully and LW is partially supported by JSPS.} }
\date{}

\maketitle

\begin{abstract}
The exceptional group $G_2$ has two maximal parabolic subgroups $P_{\mathrm{long}},\ P_{\mathrm{short}}$ corresponding to the so-called long root and short root. In this paper,
the second author introduces two zeta functions associated to $(G_2,P_{\mathrm {long}})$ and $(G_2,P_{\mathrm{short}})$ respectively, and the first author proves that these zetas satisfy the Riemann Hypothesis.
\end{abstract}
%
%%%%%%%%%%%%%%%%%%%%%%%%%%%%%%%%%%%%%%%%%%%%%%%%%%%%%%%%%%%%%%%%%%%%%%%%%%%%%%%%%%%%%%%%%%%%%%%%%%%%%%%%%%%%%%%%%%%%%%%%%%
%%
%% section 1
%%
%%%%%%%%%%%%%%%%%%%%%%%%%%%%%%%%%%%%%%%%%%%%%%%%%%%%%%%%%%%%%%%%%%%%%%%%%%%%%%%%%%%%%%%%%%%%%%%%%%%%%%%%%%%%%%%%%%%%%%%%%%
%
\section{Introduction}
Associated to a number field $F$ is the genuine high rank zeta function
$\xi_{F,r}(s)$ for every fixed $r\in\mathbb Z_{>0}$.
Being natural generalizations of (completed)
Dedekind zeta functions, these functions satisfy canonical properties for zetas as well.
Namely, they admit meromorphic continuations to the whole complex $s$-plane,
satisfy the functional equation $\xi_{F,r}(1-s)=\xi_{F,r}(s)$ and have only two singularities,
all simple poles, at $s=0,\,1$. Moreover, we expect that the Riemann Hypothesis holds for all zetas $\xi_{F,r}(s)$, namely,
all zeros of $\xi_{F,r}(s)$ lie on the central line $\mathrm {Re}(s)=\frac{1}{2}$.

Recall that $\xi_{F,r}(s)$ is defined by
\begin{equation*}
\xi_{F,r}(s):=
\Big(|\Delta_F|\Big)^{\frac{rs}{2}} \int_{\mathcal M_{F,r}}
\Big(e^{h^0(F,\Lambda)}-1\Big)\cdot\big(e^{-s}\big)^{\mathrm{deg}(\Lambda)}\,d\mu(\Lambda),
\quad \mathrm{Re}(s)>1
\end{equation*}
where $\Delta_F$ denotes the discriminant of $F$,
$\mathcal M_{F,r}$ the moduli space of semi-stable $\mathcal O_F$-lattices of rank $r$
(here $\mathcal O_F$ denotes the ring of integers),
$h^0(F,\Lambda)$ and $\mathrm{deg}(\Lambda)$ denote the 0-th geo-arithmetic cohomology
and the Arakelov degree of the lattice $\Lambda$, respectively,
and $d\mu(\Lambda)$ a certain Tamagawa type measure on $\mathcal M_{F,r}$.
Defined using high rank lattices,
these zetas then are expected to be naturally related with non-abelian aspects of number fields.
For details, see \cite{W1, W2, W3} for basic theory,
and \cite{LS} for the Riemann Hypothesis arguments.

On the other hand, algebraic groups associated to $\mathcal O_F$-lattices
are general linear group $GL$ and special linear group $SL$.
A natural question then is whether principal lattices associated to other reductive groups $G$
and their associated zeta functions can be introduced and studied.
In this paper we work with the exceptional group $G_2$. 
In contrasting  with a geo-arithmetic method used for high rank zetas (\cite{W1, W3}),
the one adopted in this paper is rather analytic (\cite{Ar1, Ar2}, \cite{JLR}).

\bigskip
For simplicity, take $F$ to be the field $\mathbb Q$ of rationals. Then, via
a Mellin transform, high rank zeta $\xi_{\mathbb Q,r}(s)$ can be written
as $$\xi_{\mathbb Q,r}(s)=\int_{\mathcal M_{\mathbb Q,r}[1]}\widehat E(\Lambda,s)\,d\mu(\Lambda),\quad \mathrm{Re}(s)>1,$$ where 
$\mathcal M_{\mathbb Q,r}[1]$ denotes the moduli space of  
$\mathbb Z$-lattices of rank $r$ and volume 1 and $\widehat E(\Lambda,s)$ the completed
 Epstein zeta functions associated to $\Lambda$.
Note that $\mathcal M_{\mathbb Q,r}[1]$ may be viewed as a compact subset in $SL(r,\mathbb Z)\backslash SL(r,\mathbb R)/SO(r)$ and 
Epstein zeta functions may be written as the relative Eisenstein series
$E^{SL(r)/P_{r-1,1}}(\bold 1;s;g)$ associated to the constant function $\bold 1$ on the maximal parabolic subgroup $P_{r-1,1}$ corresponding to the partition $r=(r-1)+1$ of $SL(r)$, we have
$$\begin{aligned}\xi_{\mathbb Q,r}(s)=&\int_{\mathcal M_{\mathbb Q,r}[1]\subset SL(r,\mathbb Z)\backslash SL(r,\mathbb R)/SO(r)}\widehat E(\Lambda,s)\,d\mu(g)\\
=&
\int_{SL(r,\mathbb Z)\backslash SL(r,\mathbb R)/SO(r)}{\bold 1}_{\mathcal M_{\mathbb Q,r}[1]}(g)\cdot\widehat E(\bold 1;s;g)\,d\mu(g)\end{aligned}$$ where ${\bold 1}_{\mathcal M_{\mathbb Q,r}[1]}(g)$ denotes the characteristic function of the compact subset $\mathcal M_{\mathbb Q,r}[1]$.

In parallel, to remedy the divergence of  integration
$$\int_{SL(r,\mathbb Z)\backslash SL(r,\mathbb R)/SO(r)}\widehat E(\bold 1;s;g)\,d\mu(g),$$ in theories of automorphic forms and trace formula, Rankin, Selberg and Arthur introduced an analytic truncation
for smooth functions $\phi(g)$ over $SL(r,\mathbb Z)\backslash SL(r,\mathbb R)/SO(r)$. Simply put, Arthur's analytic truncation
is a device to get rapidly decreasing functions from slowly increasing functions by cutting off slow growth parts near all type of cusps uniformly. Being truncations near cusps, a rather large, or better, sufficiently regular, new parameter $T$ must be introduced. In particular, when applying to Eisenstein series
$\widehat E(\bold 1;s;g)$ and to $\bold 1$ on $SL(r,\mathbb R)$, we get the truncated function $\Lambda^T\widehat E(\bold 1;s;g)$ and $(\Lambda^T\bold 1)(g)$, respectively. Consequently,
by using basic properties on Arthur's truncation (see Section 2 below),
we obtain the following well-defined integrations
$$
\begin{aligned}\int_{SL(r,\mathbb Z)\backslash SL(r,\mathbb R)/SO(r)}\Lambda^T\widehat E(\bold 1;s;g)\,d\mu(g)=&\int_{SL(r,\mathbb Z)\backslash SL(r,\mathbb R)/SO(r)}(\Lambda^T\bold 1)(g)\cdot \widehat E(\bold 1;s;g)\,d\mu(g)\\
=&\int_{\frak F(T)\subset SL(r,\mathbb Z)\backslash SL(r,\mathbb R)/SO(r)}\widehat E(\bold 1;s;g)\,d\mu(g)\end{aligned}
$$ where $\frak F(T)$ is the compact subset in (a fundamental domain of)
$SL(r,\mathbb Z)\backslash SL(r,\mathbb R)/SO(r)$ whose characteristic function is given by $(\Lambda^T\bold 1)(g)$.

\bigskip
As such, we find an analytic way to understand our high rank zetas, provided that the above analytic discussion for sufficiently positive parameter $T$ can be
further strengthened so as to work for smaller $T$, in particular, for $T=0$, as well. In general, it is very difficult (\cite{Ar1, Ar2, Ar3}). Fortunately, in the case of $SL$, this can be achieved based on an intrinsic
geo-arithmetic result, called the Micro-Global Bridge (\cite{W1, W3}), an analogue of the following basic principle in Geometric Invariant Theory for unstability: A point is not GIT stable, then there is a parabolic subgroup which destroys the stability. Consequently,
we have $$\xi_{\mathbb Q,r}(s)=\Big(\int_{G(\mathbb Z)\backslash G(\mathbb R)/K}\Lambda^T\widehat E(\bold 1;s;g)\,d\mu(g)\Big)\Big|_{T=0}.$$
This then leads to evaluation of the
special Eisenstein periods $$\int_{G(\mathbb Z)\backslash G(\mathbb R)/K}\Lambda^T\widehat E(\bold 1;s;g)\,d\mu(g),$$
and more generally the evaluation of {\it Eisenstein periods}
$$\int_{G(\mathbb Z)\backslash G(\mathbb R)/K}\Lambda^T E(\phi;\lambda;g)\,d\mu(g),$$ where 
$K$ a certain maximal compact subgroup of a reductive group $G$, $\phi$ is a $P$-level automorphic forms with $P$  parabolic,  and $E(\phi;\lambda;g)$ the relative Eisenstein 
series from $P$ to $G$ associated to $\phi$.

Unfortunately, in general, it is quite difficult to find a close formula for Eisenstein periods. But, when $\phi$ is cuspidal, then the corresponding Eisenstein period can be calculated, thanks to the work of \cite{JLR}, 
an advanced version of Rankin-Selberg \& Zagier method.

\bigskip
Back to high rank zeta functions, the bad news is that this powerful calculation cannot be applied directly, since in the specific Eisenstein series, i.e., the classical Epstein zeta, used, the function $\bold 1$, corresponding to $\phi$ in general picture, on the maximal parabolic $P_{r-1,1}$  is only $L^2$, far from being cuspidal. To overcome this technical difficulty, 
we, partially also motivated by our earlier work on the so-called abelian part of high rank zeta functions (\cite{W-2, W0}) and Venkov's trace formula for $SL(3)$ (\cite{V}), introduce Eisenstein series $E^{G/B}(\bold 1;\lambda;g)$
associated to the constant function  $\bold 1$ on $P_{1,1,\dots,1}$, the Borel, into our study, since 

1) being over the Borel, the constant function $\bold 1$ is cuspidal. So the associated Eisenstein period $\omega_{\mathbb Q}^{G;T}(\lambda)$ can be evaluated; and

2)  $E(\bold 1;s;g)$ used in high rank zetas can be realized as residues of 
$E^{G/B}(\bold 1;\lambda;g)$ along with suitable singular hyper-planes, a result essentially due to Diehl (\cite{D}). 

In fact, for 1), we have $$\omega_{\mathbb Q}^{G;T}(\lambda)=\sum_{w\in W}\Bigg(\frac{e^{\langle w\lambda-\rho,T\rangle}}{\prod_{\alpha\in\Delta_0}\langle w\lambda-\rho,\alpha^\vee\rangle}\cdot\prod_{\alpha>0,w\alpha<0}\frac{\xi_{\mathbb Q}(\langle\lambda,\alpha^\vee\rangle)}
{\xi_{\mathbb Q}(\langle\lambda,\alpha^\vee\rangle+1)}\Bigg).$$
(See Section 2 for details and unknown notations.) And for 2), we first 
know that is true for $SL(3)$ only, with the use of classical Koecher zeta functions. (See e.g., \cite{W3} for details). In believing 2) holds for general $SL(r)$, we seek the help from Henry H. Kim, among others. This proves to be quite fruitful: not only in \cite{KW}, we can offer a general formula for volume 
of truncated domain $\frak F(T)$ in the case of split, semi-simple groups, which then offers an alternative proof for Siegel-Langlands' well-known 
formula on volume of fundamental domains; 
but he brings us the paper of Diehl (\cite{D}), which deals with 
Siegel-Eisenstein series associated to the group $Sp$, from which 2) is exposed by a certain extra effort (\cite{W5}).

With all this, it is clear that there are huge difficulties in introducing and studying new
zetas associated to reductive groups $G$ geo-arithmetrically, starting from principal lattices and following the outline above 
for high rank zetas associated to $SL$. So
we decide to adapt an analytic method by focusing on the period $\omega_{\mathbb Q}^G(\lambda)$ defined by
$$\omega_{\mathbb Q}^G(\lambda):=\sum_{w\in W}\Bigg(\frac{1}{\prod_{\alpha\in\Delta_0}\langle w\lambda-\rho,\alpha^\vee\rangle}\cdot\prod_{\alpha>0,w\alpha<0}\frac{\xi_{\mathbb Q}(\langle\lambda,\alpha^\vee\rangle)}
{\xi_{\mathbb Q}(\langle\lambda,\alpha^\vee\rangle+1)}\Bigg),\ \ \mathrm{Re}\,\lambda\in\mathcal C^+.$$ Such a period, as said above, may be understood formally as
the evaluation of the Eisenstein period
$$\int_{G(F)\backslash G(\mathbb A)/K}\Lambda^TE(\bold 1;\lambda;g)d\mu(g)$$ at $T=0$ even $T$ originally is supposed to be sufficiently positive.

\bigskip
The period $\omega_{\mathbb Q}^G(\lambda)$ of $G$ over $\mathbb Q$
is of $\mathrm{rank}(G)$ variables. To get a single variable zeta out from it, totally $\mathrm{rank}(G)-1$ (linearly independent) singular hyper-planes need be chosen properly. This is done for $SL$ and $Sp$ in [W4, W5], thanks to the paper of \cite{D}. In fact, \cite{D} deals with $Sp$ only. But due to the fact that positive definite matrices are naturally associated to $\mathbb Z$-lattices and Siegel upper spaces, $SL$ can be also treated successfully with extra care. Simply put, for each $G=SL(r)$ (or $=Sp(2n)$), within the framework of classical Eisenstein series, 
there exists {\it only one} choice of  $\mathrm{rank}(G)-1$ singular hyper-planes $H_1=0,\,H_2=0,\dots,H_{\mathrm{rank}(G)-1}=0$. Moreover, after taking residues along with them, that is,  $$\mathrm{Res}_{H_1=0,\,H_2=0,\dots,H_{\mathrm{rank}(G)-1}=0}\,\omega_{\mathbb Q}^G(\lambda),$$ with suitable normalizations, we can get a new zeta $\xi_{G;\mathbb Q}(s)$ for $G$. Examples for $SL(4,\,5)$ and $Sp(4)$ shows that all these new zetas satisfy the functional equation $\xi_{G;\mathbb Q}(1-s)=\xi_{G;\mathbb Q}(s)$ and numerical tests (by MS) give supportive evidence for the RH as well. (For details, see \cite{W5}.)

At this point, the role played in new zetas $\xi_{G;\mathbb Q}(s)$ by maximal parabolic subgroups has not yet emerged. It is only after the study done for $G_2$ that we understand such a key role. Nevertheless, what we do observe from these discussions on $SL$ and $Sp$ the follows: all singular hyper-planes are taken from only a single term appeared in the period $\omega_{\mathbb Q}^G(\lambda)$, to be more precise, the term corresponding to $w=\mathrm{Id}$, the  Weyl element Identity. In other words, singular hyper-planes are taken from the denominator of the expression $$\frac{1}{\prod_{\alpha\in\Delta_0}\langle \lambda-\rho,\alpha^\vee\rangle}.$$ (Totally, there are $\mathrm{rank}(G)$ factors, among which we have carefully chosen $\mathrm{rank}(G)-1$ for $G=SL,\, Sp$.)
In particular, for the exceptional $G_2$, being a rank two group and hence an obvious  choice for our next test, this reads as
$$\frac{1}{\langle \lambda-\rho,\alpha_{\mathrm{short}}^\vee\rangle\cdot \langle \lambda-\rho,\alpha_{\mathrm{long}}^\vee\rangle}$$ where $\alpha_{\mathrm{short}}, \,\alpha_{\mathrm{long}}$ denote
the short and long roots of $G_2$ respectively. So two possibilities,

a) $\mathrm{Res}_{\langle \lambda-\rho,\alpha_{\mathrm{short}}^\vee\rangle=0}\,\omega_{\mathbb Q}^{G_2}(\lambda)$, leading to
$\xi_{\mathbb Q}^{G_2/P_{\mathrm{long}}}(s)$ after suitable normalization; and

b) $\mathrm{Res}_{\langle \lambda-\rho,\alpha_{\mathrm{long}}^\vee\rangle=0}\,\omega_{\mathbb Q}^{G_2}(\lambda)$, leading to
$\xi_{\mathbb Q}^{G_2/P_{\mathrm{short}}}(s)$ after suitable normalization.

\noindent
Here we have used the fact that there exists a natural one-to-one and onto correspondence between collection of conjugation classes of 
maximal parabolic groups and simple roots. This is the essence of  
Definition\,\&\,Proposition in Section 3, special yet very important cases of a general construction for zetas associated to reductive groups and their maximal parabolic subgroups (\cite{W5}).
 
\bigskip
As expected, similar to high rank zetas, these newly obtained zetas $\xi_{\mathbb Q}^{G_2/P}(s)$
for $G_2$ over $\mathbb Q$ prove to be canonical as well. In particular, we have the following
%
%%%%%%%%%%%%%%%%%%%%%%%%%%%%%%%%%%%%%%%%%%%%%%%%%%%%%%%%%%%%%%%%%%%%%%%%%%%%%%%%%%%%
%%
%% Theorem 1
%%
%%%%%%%%%%%%%%%%%%%%%%%%%%%%%%%%%%%%%%%%%%%%%%%%%%%%%%%%%%%%%%%%%%%%%%%%%%%%%%%%%%%%
%

\bigskip
\noindent
{\bf Theorem}
{\it Let $P=P_{\mathrm {long}}$ or $P_{\mathrm {short}}$ and $\xi_{\mathbb Q}^{G_2/P}(s)$ be the associated zeta functions.
Then}

\noindent
(1) {\it $\xi_{\mathbb Q}^{G_2/P}(s)$ are meromorphic, and admit only finite singularities, four for each, to be more precise;}

\noindent
(2) {\it $\xi_{\mathbb Q}^{G_2/P}(s)$ satisfy the standard functional equation
\begin{equation*}
\xi_{\mathbb Q}^{G_2/P}(1-s)=\xi_{\mathbb Q}^{G_2/P}(s)
\end{equation*}}

\noindent
(3) {\it All zeros of $\xi_{\mathbb Q}^{G_2/P}(s)$ lie on the central line $\mathrm{Re}\,(s)=\frac{1}{2}$.}
%\end{thm}

\bigskip
With all this said for new zetas, we next come back to point out a difference between 
high rank zetas $\xi_{{\Q},r}(s)$
and new zetas $\xi_{{\rm SL}(r);{\Q}}(s):=\xi_{\Q}^{G/P}(s)$ attached to $(G,P)=({\rm SL}(r),P_{r-1,1})$.
Roughly speaking, starting from Eisenstein series $E^{G/B}({\bold 1};\lambda;g)$,
$\xi_{{\Q},r}(s)$ corresponds to $({\rm Res} \to \int)$-ordered construction,
and new zeta functions $\xi_{{\rm SL}(r),{\Q}}(s)$
corresponds to $(\int \to {\rm Res})$-ordered construction. Here \lq\lq $({\rm Res} \to \int)$-ordered\rq\rq$\ $means that we first take the residues then take the integration, similarly,  \lq\lq $(\int \to {\rm Res})$-ordered\rq\rq$\ $means that we first take the integration then take the residues.
We have $\xi_{{\Q},2}(s)=\xi_{{\rm SL}(2),{\Q}}(s)$, since no need  taking residue.
However, in general, there is a discrepancy between $\xi_{{\Q},r}(s)$ and $\xi_{{\rm SL}(r),{\Q}}(s)$,
because of the obstruction for the exchanging of $\int$ and ${\rm Res}$. 
For example, $\xi_{{\Q},3}(s)$ has only two singularities at $s=0,1$,
but $\xi_{{\rm SL}(3),{\Q}}(s)$ has four singularities at $s=0,\,\frac{1}{3},\,\frac{2}{3},\,1$.
Simply put, new zetas $\xi_{\Q}^{GL(r)/P_{r-1,1}}(s)=\xi_{SL(r),\Q}(s)$ while close related with high rank zetas $\xi_{\mathbb Q,r}(s)$ are quite different indeed (\cite{W5}).  
Nevertheless, we expect that the distribution of zeros for $\xi_{{\rm SL}(r),{\Q}}(s)$ are quite regular as well as for $\xi_{{\Q},r}(s)$.
In fact, we have the Riemann Hypothesis for $\xi_{{\rm SL}(2),{\Q}}(s)$
(since $\xi_{{\Q},2}(s)=\xi_{{\rm SL}(2),{\Q}}(s)$), 
for $\xi_{{\rm SL}(3),{\Q}}(s)$, and for $\xi_{{\rm Sp}(4),{\Q}}(s)$ (\cite{LS, S, S2}).
All this in turn suggests that the study of new zetas $\xi_{F}^{G/P}(s)$ is not only interesting itself
but also suggestive for the study of other zetas, including Dedekind zeta functions.

\bigskip
This paper is organized as follows.
In Sections 2,\,3, we introduce various periods associated to
automorphic forms using Arthur's analytic truncations (\S2), and
define zeta functions associated to $G_2$ and its maximal parabolic subgroups (\S3).
In Sections 4,\,5, and 6, we give a proof on the corresponding Riemann Hypothesis.
%
%%%%%%%%%%%%%%%%%%%%%%%%%%%%%%%%%%%%%%%%%%%%%%%%%%%%%%%%%%%%%%%%%%%%%%%%%%%%%%%%%%%%%%%%%%%%%%%%%%%%%%%%%%%%%%%%%%%%%%%%%%
%%
%% section 2
%%
%%%%%%%%%%%%%%%%%%%%%%%%%%%%%%%%%%%%%%%%%%%%%%%%%%%%%%%%%%%%%%%%%%%%%%%%%%%%%%%%%%%%%%%%%%%%%%%%%%%%%%%%%%%%%%%%%%%%%%%%%%
%
\section{Various Periods}
In this section, we introduce various periods associated to automorphic forms using Arthur's analytic truncation.
\subsection{Automorphic Forms and Eisenstein Series}
To facilitate our ensuing discussion, we make the following preparation. For details, see e.g. \cite{MW} and/or \cite{W-1}.

Let $F$ be a number field with $\mathbb A=\mathbb A_F$ its ring of adeles.
Fix a connected reductive group $G$ defined over $F$, denote by $Z_G$ its  center.
Fix a minimal parabolic subgroup $P_0$ of $G$. Then $P_0=M_0U_0$,
where as usual we fix once and for all the Levi $M_0$ and  the unipotent radical $U_0$.
A parabolic subgroup $P$ of $G$ is called standard if $P\supset P_0$.
For such groups write $P=MU$ with $M_0\subset M$ the standard Levi and $U$ the unipotent radical.
Denote by $\mathrm{Rat}(M)$ the group of rational characters of $M$,
i.e, the morphism $M\to {\mathbb G}_m$ where ${\mathbb G}_m$ denotes the multiplicative group.
Set
\begin{equation*}
\frak a_{M,\mathbb C}^* := \mathrm{Rat}(M)\otimes_{\mathbb Z}{\mathbb C}, \qquad
\frak a_{M,\mathbb C} := \mathrm{Hom}_{\mathbb Z}(\mathrm{Rat}(M),{\mathbb C}),
\end{equation*}
and
\begin{equation*}
\frak a_M^* := \mathrm{Re}\,\frak a_M^* := \mathrm{Rat}(M)\otimes_{\mathbb Z}{\mathbb R}, \qquad
\frak a_M := \mathrm{Re}\,\frak a_M := \mathrm{Hom}_{\mathbb Z}(\mathrm{Rat}(M),{\mathbb R}).
\end{equation*}
For any $\chi \in \mathrm{Rat}(M)$,
we obtain a (real) character
$|\chi|:M({\mathbb A}) \to {\mathbb R}^*$
defined by
$m=(m_v)\mapsto m^{|\chi|}:=\prod_{v\in S} |m_v|_v^{\chi_v}$
with $|\cdot|_v$ the $v$-absolute values.
Set then
$M({\mathbb A})^1:=\cap_{\chi\in \mathrm{Rat}(M)}\mathrm{Ker}|\chi|$,
which is a normal subgroup of $M({\mathbb A})$.
Set $X_M$ be the group of complex characters which are trivial on $M({\mathbb A})^1$.
Denote by
$H_M:=\log_M:M({\mathbb A})\to \frak a_{M,\mathbb C}$
the map such that
$\forall\chi\in \mathrm{Rat}(M)\subset \frak a_{M,\mathbb C}^*,\langle\chi,
\log_M(m)\rangle:=\log(m^{|\chi|})$.
Clearly,
\begin{equation*}
M({\mathbb A})^1 = \mathrm{Ker}(\log_M); \qquad
\log_M(M({\mathbb A})/M({\mathbb A})^1) \simeq \mathrm{Re}\,\frak a_M.
\end{equation*}
Hence in particular there is a natural isomorphism
$\kappa:\frak a_{M,\mathbb C}^*\simeq X_M$.
Set
\begin{equation*}
\mathrm{Re}X_M:=\kappa (\mathrm{Re}\,\frak a_M^*),\qquad
\mathrm{Im}X_M:=\kappa (i\cdot \mathrm{Re}\,\frak a_M^*).
\end{equation*}
Moreover define our working space $X_M^G$ to be the subgroup of $X_M$
consisting of complex characters of $M({\mathbb A})/M({\mathbb A})^1$
which are trivial on $Z_{G({\mathbb A})}$.

Fix a maximal compact subgroup ${\mathbb K}$
such that for all standard parabolic subgroups $P=MU$ as above,
$P({\mathbb A})\cap{\mathbb K} = \Big(M({\mathbb A})\cap{\mathbb K}\Big) \cdot \Big(U({\mathbb A})\cap{\mathbb K}\Big)$.
Hence we get the Langlands decomposition
$G({\mathbb A}) = M({\mathbb A}) \cdot U({\mathbb A}) \cdot {\mathbb K}$.
Denote by
$m_P:G({\mathbb A}) \to M({\mathbb A})/M({\mathbb A})^1$
the map
$g=m\cdot n \cdot k\mapsto M({\mathbb A})^1 \cdot m$,
where $g \in G({\mathbb A}), m \in M({\mathbb A}), n \in U({\mathbb A})$ and $k \in {\mathbb K}$.

Fix Haar measures on $M_0({\mathbb A}), ~U_0({\mathbb A}), ~{\mathbb K}$ respectively
such that
\begin{enumerate}
\item[(1)]
the induced measure on $M(F)$ is the counting measure
and the volume of the induced measure on $M(F) \backslash M({\mathbb A})^1$ is 1.
(Recall that it is  a fundamental fact that $M(F) \backslash M({\mathbb A})^1$ is of finite volume.)

\item[(2)]
the induced measure on $U_0(F)$ is the counting measure and the volume of $U_0(F)\backslash U_0({\mathbb A})$ is 1.
(Recall that being unipotent radical, $U_0(F) \backslash U_0({\mathbb A})$ is compact.)

\item[(3)] the volume of ${\mathbb K}$ is 1.
\end{enumerate}

Such measures then also induce Haar measures via $\log_M$
to the spaces $\frak a_{M_0}, ~\frak a_{M_0}^*$, etc.
Furthermore, if we denote by $\rho_0$ the half
of the sum of the positive roots of the maximal split torus $T_0$ of the central $Z_{M_0}$ of $M_0$,
then
\begin{equation*}
f\mapsto \int_{M_0({\mathbb A})\cdot U_0({\mathbb A})\cdot
{\mathbb K}}f(mnk)\,dk\,dn\,m^{-2\rho_0}dm
\end{equation*}
defined for continuous functions with compact supports on $G({\mathbb A})$
defines a Haar measure $dg$ on $G({\mathbb A})$.
This in turn gives measures on $M({\mathbb A}), U({\mathbb A})$
and hence on $\frak a_{M}, ~\frak a_{M}^*$, $P({\mathbb A})$, etc,
for all parabolic subgroups $P$.
In particular, one checks that the following compatibility condition holds
\begin{equation*}
\int_{M_0({\mathbb A})\cdot U_0({\mathbb A})
\cdot {\mathbb K}}f(mnk)\,dk\,dn\,
m^{-2\rho_0}dm
= \int_{M({\mathbb A})\cdot U({\mathbb A})\cdot {\mathbb K}}f(mnk)\,dk\,
dn \, m^{-2\rho_P}dm
\end{equation*}
for all continuous functions $f$ with compact supports on $G({\mathbb A})$,
where $\rho_P$ denotes one half
of the sum of all positive roots of the maximal split torus $T_P$ of the central $Z_{M}$ of $M$.
For later use,
denote also by $\Delta_P$ the set of positive roots determined
by $(P,T_P)$ and $\Delta_0=\Delta_{P_0}$.

Fix an isomorphism $T_0\simeq {\mathbb G}_m^R$.
Embed ${\mathbb R}_+^*$ by the map $t\mapsto (1;t)$.
Then we obtain a natural injection
$({\mathbb R}_+^*)^R \hookrightarrow T_0({\mathbb A})$
which splits.
Denote by $A_{M_0({\mathbb A})}$ the unique connected subgroup of $T_0({\mathbb A})$
which projects onto $({\mathbb R}_+^*)^R$.
More generally, for a standard parabolic subgroup $P=MU$, set
$A_{M({\mathbb A})} := A_{M_0({\mathbb A})}\cap Z_{M({\mathbb A})}$
where as used above $Z_*$ denotes the center of the group $*$.
Clearly, $M({\mathbb A})=A_{M({\mathbb A})}\cdot M({\mathbb A})^1$.
For later use, set also
$A_{M({\mathbb A})}^G := \{a\in A_{M({\mathbb A})}:\log_Ga=0\}$.
Then $A_{M({\mathbb A})}=A_{G({\mathbb A})} \oplus A_{M({\mathbb A})}^G$.

Note that ${\mathbb K}$ and $U(F) \backslash U({\mathbb A})$ are all compact,
and $M(F) \backslash M({\mathbb A})^1$ is of finite volume.
With the Langlands decomposition
$G({\mathbb A}) = U({\mathbb A}) M({\mathbb A}) {\mathbb K}$ in mind,
the reduction theory for $G(F) \backslash G({\mathbb A})$
or more generally for $P(F) \backslash G({\mathbb A})$
is reduced to that for $A_{M({\mathbb A})}$
since $Z_G(F) \cap Z_{G({\mathbb A})} \backslash Z_{G({\mathbb A})}\cap G({\mathbb A})^1$
is compact as well.
As such, for $t_0\in M_0({\mathbb A})$ set
\begin{equation*}
A_{M_0({\mathbb A})}(t_0)
:= \{ a \in A_{M_0({\mathbb A})} : a^\alpha>t_0^\alpha,  ~\forall \alpha \in \Delta_0 \}.
\end{equation*}
Then, for a fixed compact subset $\omega\subset P_0({\mathbb A})$,
we have the corresponding Siegel set
\begin{equation*}
S(\omega;t_0) := \{ p \cdot a \cdot k : p \in \omega, ~a \in A_{M_0({\mathbb A})} (t_0),~k \in {\mathbb K} \}.
\end{equation*}
In particular, the classical reduction theory may be restated as,
for big enough $\omega$ and small enough $t_0$,
i.e, $t_0^\alpha$ is very close to 0 for all $\alpha \in \Delta_0$,
$G({\mathbb A})=G(F)\cdot S(\omega;t_0)$.
More generally set
\begin{equation*}
A_{M_0({\mathbb A})}^P(t_0)
:= \{ a \in A_{M_0({\mathbb A})} : a^\alpha > t_0^\alpha, ~\forall \alpha \in \Delta_0^P\},
\end{equation*}
and
\begin{equation*}
S^P(\omega;t_0)
:= \{ p \cdot a \cdot k : p \in \omega, ~a \in A_{M_0({\mathbb A})}^P(t_0), ~k \in {\mathbb K} \}.
\end{equation*}
Then similarly as above
for big enough $\omega$ and small enough $t_0$,
$G({\mathbb A})=P(F) \cdot S^P(\omega;t_0)$.
(Here $\Delta_0^P$ denotes the set of positive roots for $(P_0 \cap M,T_0)$.)

Fix an embedding $i_G : G \hookrightarrow SL_n$ sending $g$ to $(g_{ij})$.
Introducing a height function on $G({\mathbb A})$ by setting
$\|g\| := \prod_{v\in S} \sup \{|g_{ij}|_v : \forall i,j\}$.
It is well known that up to $O(1)$, height functions are unique.
This implies that the following growth conditions do not depend on the height function we choose.

A function $f:G({\mathbb A}) \to {\mathbb C}$ is said to have moderate growth
if there exist $c,r\in {\mathbb R}$ such that $|f(g)| \leq c \cdot \| g \|^r$ for all $g \in G({\mathbb A})$.
Similarly, for a standard parabolic subgroup $P=MU$,
a function $f:U({\mathbb A})M(F) \backslash G({\mathbb A}) \to {\mathbb C}$
is said to have moderate growth if there exist
$c,r \in {\mathbb R},\lambda \in \mathrm{Re}X_{M_0}$
such that for any
$a \in A_{M({\mathbb A})},k\in {\mathbb K}, ~m \in M({\mathbb A})^1\cap S^P(\omega;t_0)$,
\begin{equation*}
|f(amk)| \leq c \cdot \|a\|^r \cdot m_{P_0}(m)^\lambda.
\end{equation*}

By contrast, a function $f: S(\omega;t_0)\to {\mathbb C}$ is said to be rapidly decreasing
if there exists $r>0$ and for all $\lambda\in \mathrm{Re}X_{M_0}$
there exists $c>0$ such that for
$a \in A_{M({\mathbb A})}, ~g \in G({\mathbb A})^1\cap S(\omega;t_0)$,
$|\phi(ag)| \leq c \cdot \|a\|\cdot m_{P_0}(g)^\lambda$.
And a function $f:G(F) \backslash G({\mathbb A}) \to {\mathbb C}$
is said to be rapidly decreasing if $f|_{S(\omega;t_0)}$ is so.

Also a function $f:G({\mathbb A}) \to {\mathbb C}$ is said to be smooth
if for any $g=g_f \cdot g_\infty \in G({\mathbb A}_f) \times G({\mathbb A}_\infty)$,
there exist open neighborhoods $V_*$ of $g_*$ in $G({\mathbb A})$
and a $C^\infty$-function $f':V_\infty\to {\mathbb C}$
such that $f(g_f'\cdot g_\infty')=f'(g_\infty')$ for all $g_f'\in V_f$ and $g_\infty'\in V_\infty$.

By definition, a function
$\phi: U({\mathbb A})M(F) \backslash G({\mathbb A}) \to {\mathbb C}$
is called {\it automorphic}
if

\begin{enumerate}
\item[(i)] $\phi$ has moderate growth;

\item[(ii)] $\phi$ is smooth;

\item[(iii)] $\phi$ is ${\mathbb K}$-finite,
i.e, the ${\mathbb C}$-span of all $\phi(k_1\cdot *\cdot k_2)$
parametrized by $(k_1,k_2)\in {\mathbb K}\times {\mathbb K}$
is finite dimensional; and

\item[(iv)] $\phi$ is $\frak z$-finite,
i.e, the ${\mathbb C}$-span of all $\delta(X)\phi$
parametrized by all $X \in \frak z$ is finite dimensional.
Here $\frak z$ denotes the center of the universal enveloping algebra
$\frak u:=\frak U(\text{Lie}G({\mathbb A}_\infty))$ of the Lie algebra of
$G({\mathbb A}_\infty)$ and $\delta(X)$ denotes the derivative of $\phi$ along $X$.
\end{enumerate}

\noindent
Set $A(U({\mathbb A})M(F)\backslash G({\mathbb A}))$ be the space of automorphic
forms on $U({\mathbb A})M(F)\backslash G({\mathbb A})$.

For a measurable locally $L^1$-function
$f:U(F) \backslash G({\mathbb A}) \to {\mathbb C}$,
define its {\it constant term} along with the standard parabolic subgroup $P=UM$ to be
$f_P:U({\mathbb A}) \backslash G({\mathbb A})\to {\mathbb C}$
given by
$g \mapsto \int_{U(F) \backslash G({\mathbb A})}f(ng)dn$.
Then an automorphic form
$\phi \in A(U({\mathbb A})M(F) \backslash G({\mathbb A}))$
is called a {\it cusp form}
if for any standard parabolic subgroup $P'$ properly contained in $P$, $\phi_{P'}\equiv 0$.
Denote by $A_0(U({\mathbb A})M(F) \backslash G({\mathbb A}))$
the space of cusp forms on $U({\mathbb A})M(F) \backslash G({\mathbb A})$.
One checks easily that

\begin{enumerate}
\item[(i)] all cusp forms are rapidly decreasing; and hence

\noindent
\item[(ii)] there is a natural pairing
\begin{equation*}
\langle \cdot,\cdot \rangle :
A_0(U({\mathbb A})M(F) \backslash G({\mathbb A})) \times A(U({\mathbb A})M(F) \backslash G({\mathbb A})) \to {\mathbb C}
\end{equation*}
defined by
$ \langle \psi,\phi \rangle :=
\int_{Z_{M({\mathbb A})}U({\mathbb A})M(F) \backslash G({\mathbb A})}
\psi(g)\bar\phi(g)\,dg$.
\end{enumerate}

For an  automorphic form
$\phi\in A(U({\mathbb A})M(F)\backslash G({\mathbb A}))$,
define the associated {\it Eisenstein series}
$E(\phi,\lambda):G(F)\backslash G({\mathbb A})\to {\mathbb C}$
by
\begin{equation*}
E(\phi,\lambda)(g):=\sum_{\delta\in P(F)\backslash G(F)}
\phi(\delta g)\cdot m_P(\delta g)^{\lambda+\rho_P}.
\end{equation*}
Then one checks that there is an open cone
${\mathcal C}\subset \mathrm{Re}X_M^G$
such that if $\mathrm{Re}\lambda \in {\mathcal C}$,
$E(\phi,\lambda)(g)$ converges uniformly for $g$
in a compact subset of $G({\mathbb A})$ and $\lambda$ in an open neighborhood of 0 in $X_M^G$.
For example, if $\phi$ is cuspidal,
we may even take ${\mathcal C}$ to be the cone
$\{ \lambda \in \mathrm{Re} X_M^G : \langle \lambda, ~\alpha^\vee\rangle>0,~\forall\alpha\in \Delta_P^G\}$.
As a direct consequence,
then $E(\phi,\lambda)\in A(G(F)\backslash G({\mathbb A}))$.
That is, it is an automorphic form.

We end this discussion by introducing intertwining operators.
For $w \in W$ the Weyl group of $G$, fix once and for all representative $w \in G(F)$ of $w$.
Set $M':=wMw^{-1}$ and denote the associated parabolic subgroup by $P'=U'M'$.
As usual, define the associated intertwining operator $M(w,\lambda)$ by
\begin{equation*}
\begin{aligned}
&\Big(M(w,\lambda)\phi\Big)(g)
:= m_{P'}(g)^{w\lambda+\rho_{P'}} \\
&\times \int_{U'(F)\cap wU(F)w^{-1} \backslash U'({\mathbb  A})}
\phi(w^{-1}n'g)\cdot m_P(w^{-1}n'g)^{\lambda+\rho_P}dn', \quad
\forall g\in G({\mathbb  A}).
\end{aligned}
\end{equation*}

\subsection{Arthur's Analytic Truncation}

Let $P$ be a (standard) parabolic subgroup of $G$.
Write $T_P$ for the maximal split torus in the center of $M_P$
and $T_P'$ for the maximal quotient split torus of $M_P$.
Set $\tilde{\frak a}_P:=X_*(T_P)\otimes \mathbb R$
and denote its real dimension by $d(P)$,
where $X_*(T)$ is the lattice of 1-parameter subgroups in the torus $T$.
Then it is known that $\tilde{\frak a}_P=X_*(T_P') \otimes \mathbb R$ as well.
The two descriptions of $\tilde {\frak a}_P$ show that
if $Q\subset P$ is a parabolic subgroup,
then there is a canonical injection
$\tilde {\frak a}_{P} \hookrightarrow \tilde{\frak a}_{Q}$
and a natural surjection
$\tilde {\frak a}_{Q} \twoheadrightarrow \tilde{\frak a}_{P}$.
We thus obtain a canonical decomposition
$\tilde {\frak a}_{Q}=\tilde{\frak a}_Q^{P}\oplus \tilde {\frak a}_{P}$
for a certain subspace $\tilde{\frak a}_Q^{P}$ of $\tilde {\frak a}_{Q}$.
In particular, $\tilde{\frak a}_G$ is a summand of
$\tilde {\frak a}=\tilde{\frak a}_P$ for all $P$.
Set
$\frak a_{P}:= \tilde{\frak a}_{P}/\tilde {\frak a}_{G}$
and
$\frak a_{Q}^P:= \tilde{\frak a}_Q^{P}/ \tilde{\frak a}_{G}$.
Then we have
\begin{equation*}
\frak a_{Q}=\frak a_Q^{P}\oplus\frak a_{P}
\end{equation*}
and $\frak a_P$ is canonically identified as a subspace of $\frak a_Q$.
Set $\frak a_0 := \frak a_{P_0}$ and $\frak a_0^P = \frak a_{P_0}^P$
then we also have $\frak a_0=\frak a_0^P\oplus \frak a_P$ for all $P$.

Dually we have spaces
$\frak a_0^*, ~\frak a_P^*, ~\bigl(\frak a_0^P \bigr)^*$,
(where for a real space $V$, write $V^*$ its dual space over $\mathbb R$,)
and hence the decompositions
$\frak a_0^* = \bigl( \frak a_0^Q \bigr)^* \oplus \bigl(\frak a_Q^P \bigr)^* \oplus \frak a_P^*.$

So $\frak a_P^*=X(M_P)\otimes \mathbb R$ with $X(M_P)$ the group
$\mathrm {Hom}_F \bigl( M_P,GL(1) \bigr)$
i.e., collection of characters on $M_P$.
It is known that
$\frak a_P^*=X(A_P)\otimes \mathbb R$
where $A_P$ denotes the split component of the center of $M_P$.
Clearly, if $Q\subset P$, then $M_Q\subset M_P$ while $A_P\subset A_Q$.
Thus via restriction, the above two expressions of $\frak a_P^*$
also naturally induce an injection $\frak a_P^*\hookrightarrow \frak a_Q^*$
and a surjection $\frak a_Q^* \twoheadrightarrow \frak a_P^*$,
compatible with the decomposition
$\frak a_Q^* = \bigl( \frak a_Q^P \bigr)^* \oplus \frak a_P^*$.

As usual, let $\Delta_0$ and $\widehat\Delta_0$
be the subsets of simple roots and simple weights in $\frak a_0^*$ respectively.
Write $\Delta_0^\vee$ (resp. $\widehat\Delta^\vee_0$)
for the basis of $\frak a_0$ dual to $\widehat\Delta_0$ (resp. $\Delta_0$).
Being the dual of the collection of simple weights (resp. of simple roots),
$\Delta_0^\vee$ (resp. $\widehat\Delta_0^\vee$)
is the set of coroots (resp. coweights).

For every $P$, let $\Delta_P \subset \frak a_0^*$ be the set of non-trivial
{\it restrictions} of elements of $\Delta_0$ to $\frak a_P$.
Denote the dual basis of $\Delta_P$ by $\widehat\Delta_P^\vee$.
For each $\alpha \in \Delta_P$,
let $\alpha^\vee$ be the projection of $\beta^\vee$ to $\frak a_P$,
where $\beta$ is the root in $\Delta_0$ whose restriction to $\frak a_P$ is $\alpha$.
Set $\Delta_P^\vee := \bigl\{ \alpha^\vee:\alpha\in\Delta_P \bigr\}$,
and define the dual basis of $\Delta_P^\vee$ by $\widehat\Delta_P$.

More generally, if $Q\subset P$,
write $\Delta_Q^P$ to denote the {\it subset} $\alpha\in\Delta_Q$
appearing in the action of $T_Q$ in the unipotent radical of $Q\cap M_P$.
(Indeed, $M_P\cap Q$ is a parabolic subgroup of $M_P$ with nilpotent radical $N_Q^P := N_Q\cap M_P$.
Thus $\Delta_Q^P$ is simply the set of roots of the parabolic subgroup $(M_P \cap Q,A_Q)$.
And one checks that the map $P \mapsto \Delta_Q^P$ gives a natural bijection
between parabolic subgroups $P$ containing $Q$ and subsets of $\Delta_Q$.)
Then $\frak a_P$ is the subspace of $\frak a_Q$ annihilated by $\Delta_Q^P$.
Denote by $(\widehat\Delta^\vee)_Q^P$ the dual of $\Delta_Q^P$.
Let $(\Delta_Q^P)^\vee := \bigl\{ \alpha^\vee : \alpha \in \Delta_Q^P \bigr\}$
and denote by $\widehat\Delta_Q^P$ the dual of $(\Delta_Q^P)^\vee$.

Also we extend the linear functionals in $\Delta_Q^P$ and $\widehat\Delta_Q^P$
to elements of the dual space $\frak a_0^*$
by means of the canonical projection from $\frak a_0$ to $\frak a_Q^P$
given by the decomposition
$\frak a_0 = \frak a_0^Q \oplus \frak a_Q^P \oplus \frak a_P$.
Let $\widehat\tau_Q^P$ be the characteristic function of the {\it positive cone}
\begin{equation*}
\begin{aligned}
\Bigl\{ H \in \frak a_0: \langle \varpi, H \rangle & >0, ~\forall \varpi \in \widehat\Delta_Q^P \Bigr\}\\
= &\, \frak a_0^Q\oplus
\Big\{ H \in \frak a_Q^P: \langle \varpi, H \rangle >0
~\text{for all}\ \varpi\in\widehat\Delta_Q^P\Big\}
\oplus\frak a_P.
\end{aligned}
\end{equation*}
Denote  $\widehat\tau_P^G$ simply by $\widehat\tau_P$.

Recall that an element $T\in\frak a_0$ is called {\it sufficiently regular},
if $\alpha(T)\gg 0$ for any $\alpha \in \Delta_0$.
Fix then a suitably regular point $T \in \frak a_0$.
If $\phi$ is a continuous function on $G(F) \backslash G(\mathbb A)^1$,
define {\it Arthur's analytic truncation}
$\bigl( \Lambda^T\phi \bigr)(x)$ to be the function
\begin{equation*}
\bigl( \Lambda^T\phi \bigr)(x):=
\sum_P(-1)^{\mathrm{dim} (A/Z)}\sum_{\delta\in P(F)\backslash G(F)}
\phi_P(\delta x)\cdot\hat\tau_P \bigl( H(\delta x)-T \bigr),
\end{equation*}
where
\begin{equation*}
\phi_P(x):=\int_{N(F) \backslash N(\mathbb A)}\phi(nx)\,dn
\end{equation*}
denotes the constant term of $\phi$ along $P$,
and the sum is over all (standard) parabolic subgroups.

Note that all parabolic subgroups of $G$ can be obtained from standard parabolic subgroups
by taking conjugations with elements from $P(F)\backslash G(F)$.
So we have:
\vskip 0.30cm
\begin{enumerate}
\item[(a)]
$\displaystyle{\bigl( \Lambda^T\phi \bigr)(x)=\sum_P(-1)^{\mathrm{dim} (A/Z)}
\phi_P(x) \cdot\hat\tau_P \bigl( H(x)-T \bigr),}$
{\it where the sum is over all, both standard and non-standard, parabolic subgroups};

\item[(b)] {\it If $\phi$ is a cusp form, then $\Lambda^T\phi=\phi$}.
\end{enumerate}

Fundamental properties of Arthur's analytic truncation may be summarized as follows:
%
%%%%%%%%%%%%%%%%%%%%%%%%%%%%%%%%%%%%%%%%%%%%%%%%%%%%%%%%%%%%%%%%%%%%%%%%%%%%%%%%%%%%
%%
%% Theorem 3, Arthur
%%
%%%%%%%%%%%%%%%%%%%%%%%%%%%%%%%%%%%%%%%%%%%%%%%%%%%%%%%%%%%%%%%%%%%%%%%%%%%%%%%%%%%%
%
\begin{thm}[Arthur~\cite{Ar1,Ar2}] \label{thm_22}
For sufficiently regular $T$ in $\frak a_0$,

\noindent
{\rm (1)} Let $\phi:G(F)\backslash G(\mathbb A)\to\mathbb C$ be a locally $L^1$-function.
Then
\begin{equation*}
\Lambda^T\Lambda^T\phi(g)=\Lambda^T\phi(g)
\end{equation*}
for almost all $g$.
If $\phi$ is also locally bounded,
then the above is true for all $g$;

\noindent
{\rm (2)} Let $\phi_1,\,\phi_2$ be two locally $L^1$ functions on $G(F)\backslash G(\mathbb A)$.
Suppose that $\phi_1$ is of moderate growth and $\phi_2$ is  rapidly decreasing.
Then
\begin{equation*}
\int_{Z_{G(\mathbb A)} G(F) \backslash G(\mathbb A)}
\overline{\Lambda^T \phi_1(g)} \cdot \phi_2(g) \,dg
= \int_{Z_{G(\mathbb A)}G(F)\backslash G(\mathbb A)}
\overline{\phi_1(g)} \cdot \Lambda^T \phi_2(g) \, dg;
\end{equation*}

\noindent
{\rm (3)} Let $K_f$ be an open compact subgroup of $G(\mathbb A_f)$,
and $r, ~r'$ are two positive real numbers.
Then there exists a finite subset
$\bigl\{ X_i:i=1,2,\ldots,N \bigr\}\subset\mathcal U$,
the universal enveloping algebra of $\frak g_\infty$,
such that the following is satisfied:
Let $\phi$ be a smooth function on $G(F)\backslash G(\mathbb A)$,
right invariant under $K_f$ and let $a\in A_{G(\mathbb A)},\ g \in G(\mathbb A)^1 \cap S$.
Then
\begin{equation*}
\Big|\Lambda^T\phi(ag)\Big|\leq \|g\|^{-r}\sum_{i=1}^N
\sup\Big\{|\delta(X_i)\phi(ag')|\,\|g'\|^{-r'}:g'\in G(\mathbb A)^1\Big\},
\end{equation*}
where $S$ is a  Siegel domain with respect to $G(F)\backslash G(\mathbb A)$.
\end{thm}

\subsection{Arthur's Periods}
\vskip 0.30cm
Fix a sufficiently regular $T \in \frak a_0$ and let $\phi$ be an automorphic form of $G$.
Then, $\Lambda^T\phi$ is rapidly decreasing, and hence integrable.
In particular, the integration
\begin{equation*}
A(\phi;T):=
\int_{G(F) \backslash G(\mathbb A)} \Lambda^T \phi(g) \,dg
\end{equation*}
makes sense.
We claim that $A(\phi;T)$ can be written as an integration of
the original automorphic form $\phi$ over a certain compact subset.

To start with, note that for Arthur's analytic truncation $\Lambda^T$,
we have $\Lambda^T\circ \Lambda^T=\Lambda^T$.
Hence,
\begin{equation*}
\begin{aligned}
A(\phi;T)
= & \int_{Z_{G({\mathbb A})}G(F) \backslash G({\mathbb A})} \Lambda^T \phi\ d\mu(g)\\
= & \int_{Z_{G({\mathbb A})}G(F) \backslash G({\mathbb A})} \Lambda^T \bigl( \Lambda^T \phi \bigr)(g)\ d\mu(g).
\end{aligned}
\end{equation*}
Moreover, by the self-adjoint property,
for the constant function $\bold 1$ on $G(\mathbb A)$,
\begin{equation*}
\begin{aligned}
& \int_{Z_{G({\mathbb A})}G(F) \backslash G({\mathbb A})}
{\bf 1}(g) \cdot \Lambda^T \bigl( \Lambda^T \phi \bigr)(g)\ d\mu(g) \\
& \quad =  \int_{Z_{G({\mathbb A})}G(F)\backslash G({\mathbb A})}
\bigl( \Lambda^T {\bf 1} \bigr)(g) \cdot
\bigl( \Lambda^T \phi \bigr)(g)\ d\mu(g) \\
& \quad =  \int_{Z_{G({\mathbb A})}G(F)\backslash G({\mathbb A})}
\Lambda^T \bigl( \Lambda^T {\bf 1} \bigr)(g) \cdot \phi(g)\ d\mu(g),
\end{aligned}
\end{equation*}
since $\Lambda^T\phi$ and $\Lambda^T{\bf 1}$ are rapidly decreasing.
Therefore, using $\Lambda^T\circ \Lambda^T=\Lambda^T$ again,
we arrive at
\begin{equation*}
A(\phi;T)
= \int_{Z_{G({\mathbb A})}G(F)\backslash G({\mathbb A})}
\Lambda^T{\bold 1}(g)\cdot \phi(g)\ d\mu(g).\eqno(*)
\end{equation*}

To go further, let us give a much more detailed study of Arthur's analytic
truncation for the constant function ${\bold 1}$.
Introduce the truncated subset $\Sigma(T)$
of the space $Z_{G({\mathbb A})}G(F)\backslash G({\mathbb A})$ by
\begin{equation*}
\Sigma(T) :=
\Big\{ g\in Z_{G({\mathbb A})}G(F)\backslash G({\mathbb A}): \Lambda^T{\bf 1}(g)=1 \Big\}.
\end{equation*}
%
%%%%%%%%%%%%%%%%%%%%%%%%%%%%%%%%%%%%%%%%%%%%%%%%%%%%%%%%%%%%%%%%%%%%%%%%%%%%%%%%%%%%
%%
%% Proposition 1, Arthur
%%
%%%%%%%%%%%%%%%%%%%%%%%%%%%%%%%%%%%%%%%%%%%%%%%%%%%%%%%%%%%%%%%%%%%%%%%%%%%%%%%%%%%%
%
\begin{prop}[Arthur~\cite{Ar3}] \label{prop_22_1}
For sufficiently regular $T\in\frak a_0$,
$\Lambda^T{\bold 1}$ is the characteristic function
of a compact subset of $Z_{G({\mathbb A})}G(F)\backslash G({\mathbb A})$.
In particular, $\Sigma(T)$ is compact.
\end{prop}

Consequently,
$\qquad \int_{Z_{G({\mathbb A})}G(F)\backslash G({\mathbb A})}\Lambda^T\phi(g)\,
d\mu(g)=$
\begin{equation*}
\int_{Z_{G({\mathbb A})}G(F)\backslash G({\mathbb A})}
\Lambda^T{\bold 1}(g)\cdot\phi(g)\, d\mu(g)=\int_{\Sigma(T)}\phi(g)\, d\mu(g).
\end{equation*}
That is to say, we have obtained the following:
%
%%%%%%%%%%%%%%%%%%%%%%%%%%%%%%%%%%%%%%%%%%%%%%%%%%%%%%%%%%%%%%%%%%%%%%%%%%%%%%%%%%%%
%%
%% Proposition 2
%%
%%%%%%%%%%%%%%%%%%%%%%%%%%%%%%%%%%%%%%%%%%%%%%%%%%%%%%%%%%%%%%%%%%%%%%%%%%%%%%%%%%%%
%
\begin{prop} \label{prop_22_2}
For a sufficiently regular $T\in\frak a_0$
and an automorphic form $\phi$ on $G(F)\backslash G(\mathbb A)$,
\begin{equation*}
\int_{\Sigma(T)}\phi(g)\,
d\mu(g)=\int_{Z_{G({\mathbb A})}G(F)\backslash G({\mathbb A})}\Lambda^T\phi(g)\, d\mu(g).
\end{equation*}
\end{prop}
\noindent
It is because of this result that we call $\int_{G(F)\backslash G({\mathbb A})^1}\Lambda^T\phi(g)\, d\mu(g)$
the \emph{Arthur period} for $\phi$.

\subsection{Eisenstein Periods}

Let $P$ be a (standard) parabolic subgroup of $G$ with Levi decomposition $P=MU$
and $\phi \in A(U(\mathbb A)M(F)\backslash G(\mathbb A))$ an $M$-level automorphic form.
Then the associated Eisenstein series
$E(\phi;\lambda)(g) := \sum_{\delta\in P(F)\backslash G(F)}
\phi(\delta g)\cdot m_P(\delta g)^{\lambda+\rho_P} \in A(G(F)\backslash G(\mathbb A)$
is a $G$-level automorphic form.
Thus for a sufficiently positive $T\in \frak a_0$,
we obtain a well-defined Arthur period
\begin{equation*}
\int_{Z_{G({\mathbb A})}G(F)\backslash G(\mathbb A)}\wedge^TE(\phi;\lambda)(g)\,d\mu(g).
\end{equation*}
Due to the obvious importance, we call such an Arthur period an \emph{Eisenstein period}.

In general, Eisenstein periods are quite difficult to be evaluated.
However, if $\phi$ is cuspidal, we have the following result of \cite{JLR},
an advanced version of the Rankin-Selberg \& Zagier method.
%
%%%%%%%%%%%%%%%%%%%%%%%%%%%%%%%%%%%%%%%%%%%%%%%%%%%%%%%%%%%%%%%%%%%%%%%%%%%%%%%%%%%%
%%
%% Theorem 4
%%
%%%%%%%%%%%%%%%%%%%%%%%%%%%%%%%%%%%%%%%%%%%%%%%%%%%%%%%%%%%%%%%%%%%%%%%%%%%%%%%%%%%%
%
\begin{thm}[\cite{JLR}]
Fix a sufficiently positive $T \in \frak a_0^+$.
Let $P=MU$ be a parabolic subgroup and $\phi$ a $P$-level cusp form.
Then the Eisenstein period
$\int_{G(F)\backslash G(\mathbb A)}\Lambda^TE(\lambda,\phi)(g)\,dg$
is equal to

\noindent
{\rm (1)} $0$ if $P\not=P_0$ is not minimal; {\rm and}

\noindent
{\rm (2)} $\displaystyle{\mathrm{Vol}\Big(\Big\{\sum_{\alpha\in\Delta_0}
a_\alpha\alpha^\vee:a_\alpha\in[0,1)\Big\}\Big)}$

\begin{equation*}
\times{\displaystyle\sum_{w\in W}}\frac{e^{\langle w\lambda-\rho,T\rangle}}
{\prod_{\alpha\in\Delta_0}\langle w\lambda-\rho,\alpha^\vee\rangle}
\cdot\int_{M_0(F)\backslash M_0(\mathbb A)^1\times K}\Big(M(w,\lambda)
\phi\Big)(mk)\,dm\,dk,\
\end{equation*}
if $P=P_0=M_0U_0$ is minimal.
\end{thm}

\subsection{Periods for $G$ over $F$}

Now we focus on the expression
\begin{equation*}
\sum_{w\in W} \frac{e^{\langle w\lambda-\rho,T\rangle}}
{\prod_{\alpha\in\Delta_0}\langle w\lambda-\rho,\alpha^\vee\rangle}
\times\int_{M_0(F)\backslash M_0(\mathbb A)^1\times K}\Big(M(w,\lambda)
\phi\Big)(mk)\,dm\,dk,\eqno(*)
\end{equation*}
for a cusp form $\phi$ at the level of the Borel.
\emph{Motivated by our study of high rank zetas} (\cite{W1,W3,W4,W5}),
we make the following two simplifications:

\begin{enumerate}
\item[(1)]
Take $T=0$. Recall that in the discussion so far,
$T$ has been assumed to be sufficiently positive.
However, $(*)$ makes sense even when $T=0$; and

\item[(2)]
Take $\phi\equiv {\bf 1}$, the constant function one on the Borel.
Recall that in general for a standard $P=MU$,
the constant function ${\bf 1}$ is only $L^2$ on $M$.
But for the Borel, ${\bf 1}$ is cuspidal.
\end{enumerate}

With all these preparations, we are ready to introduce our first main definition.
\vskip 0.30cm
\noindent
\begin{defi}
The \emph{period $\omega_{F}^G(\lambda)$ of G over $F$} is defined by
\begin{equation*}
\omega_{F}^G(\lambda):={\displaystyle\sum_{w\in W}}
\Bigg(\frac{1}
{\prod_{\alpha\in\Delta_0}\langle w\lambda-\rho,\alpha^\vee\rangle}\times M(w,\lambda)\Bigg),
\end{equation*}
where
$M(w,\lambda)$ denotes the quantity
\begin{equation*}
m_{P'}(e)^{w\lambda+\rho_{P'}}\cdot\int_{U'(F)\cap wU(F)w^{-1}\backslash U'({\mathbb  A})}
 m_P(w^{-1}n')^{\lambda+\rho_P}dn'
\end{equation*}
where
$M':=wMw^{-1}$ and $P'=U'M'$ denotes the associated parabolic subgroup.
\end{defi}

In particular, for $G=G_2$,  by the Gindikin-Karpelevich formula (\cite{L2}),
we have
\begin{equation*}
M(w,\lambda) = \prod_{\alpha>0, \, w\alpha<0}
\frac {\xi\big(\langle\lambda, \alpha^{\vee}\rangle\big)}{\xi\big(\langle \lambda, \alpha^{\vee}\rangle+1\big)}.
\end{equation*}
Here $\xi(s):=\pi^{-\frac s2}\Gamma(\frac s2)\zeta(s)$ with $\zeta(s)$ the Riemann zeta function.
Consequently,
\begin{equation*}
\omega^{G_2}_{\mathbb Q}(\lambda) := \sum_{w \in W}
\Bigg(
\frac{1}{\prod_{\alpha\in\Delta_0}\langle w\lambda-\rho,\alpha^\vee\rangle} \times \prod_{\alpha>0, \, w\alpha<0}
\frac {\xi\big(\langle\lambda, \alpha^{\vee}\rangle\big)}{\xi\big(\langle \lambda, \alpha^{\vee}\rangle+1\big)}\Bigg).
\eqno(**)
\end{equation*}
%
%%%%%%%%%%%%%%%%%%%%%%%%%%%%%%%%%%%%%%%%%%%%%%%%%%%%%%%%%%%%%%%%%%%%%%%%%%%%%%%%%%%%%%%%%%%%%%%%%%%%%%%%%%%%%%%%%%%%%%%%%%
%%
%% section 3
%%
%%%%%%%%%%%%%%%%%%%%%%%%%%%%%%%%%%%%%%%%%%%%%%%%%%%%%%%%%%%%%%%%%%%%%%%%%%%%%%%%%%%%%%%%%%%%%%%%%%%%%%%%%%%%%%%%%%%%%%%%%%
%
\section{Zetas for $G_2$}
In this section, we introduce zeta functions associated to $(G_2, P_{\mathrm{short}})$ and $(G_2, P_{\mathrm{long}})$ using the period of $G_2$ introduced in Section 2.
\subsection{Period for $G_2$ over $\mathbb Q$}

Let $G$ be the exceptional group $G_2$.
It is simply connected and adjoint.
Fix a maximal split torus $T$ in $G$ and a Borel subgroup $B$ containing $T$.
Then we obtain two simple roots,
the short root  $\alpha$ and the long root $\beta$. So $\Delta_0=\{\alpha,\beta\}$ and
and all positive roots are given by
\begin{equation*}
\Phi^+=\{\alpha,\beta,\alpha+\beta,2\alpha+\beta, 3\alpha+\beta, 3\alpha+2\beta\}.
\end{equation*}
Denote by $P_{\mathrm{long}}=P_{\beta}=P_1$ and $P_{\mathrm{short}}=P_\alpha=P_2$ the maximal standard parabolic subgroups
attached to $\Delta_0\backslash\{\beta\}$ and $\Delta_0\backslash\{\alpha\}$, respectively. (See e.g., \cite{H})

Choose a parametrization $t:\mathbb Q^*\times\mathbb Q^*\to T, (a,b)\mapsto t(a,b)$
defined by $\alpha(t(a,b))=ab^{-1}, ~\beta(t(a,b))=a^{-1}b^2$.
Then the actions of remaining positive roots are given by
\begin{equation*}
\begin{aligned}
&(\alpha+\beta)(t(a,b)) = b,\qquad (2\alpha+\beta)(t(a,b)) = a, \\
&(3\alpha+\beta)(t(a,b)) = a^2b^{-1},\qquad (3\alpha+2\beta)(t(a,b)) = ab
\end{aligned}
\end{equation*}
and the corresponding coroot are given by
\begin{equation*}
\begin{aligned}
&\alpha^\vee(x)=t(x,x^{-1}),\qquad
\beta^\vee(x)=t(1,x),\qquad (\alpha+\beta)^\vee(x)=t(x,x^2),\\
&(2\alpha+\beta)^\vee(x)=t(x^2,x),\quad (3\alpha+\beta)^\vee(x)=t(x,1),\quad (3\alpha+2\beta)^\vee(x)=t(x,x).
\end{aligned}
\end{equation*}

Let $X(T)$ be the character group of $T$ and $\frak a_{\mathbb C}^* = X(T) \otimes \mathbb C$ its complexification.
We introduce coordinates in $\frak a_{\mathbb C}^*$ with respect to the basis $2\alpha+\beta,~\alpha+\beta$.
Thus point $(z_1,z_2)\in\mathbb C^2$ corresponds to the character $\lambda = z_1(2\alpha+\beta) + z_2(\alpha+\beta)$.
(The coordinate is chosen for the reason to make $\lambda(t(a,b)) = |a|^{z_1}|b|^{z_2}$ take the simplest form.)
As such, then $\rho:=\rho_B:= 5\alpha + 3\beta$
and $\mathcal C^+$ of the positive Weyl chamber in $\frak a_{\mathbb C}^*$ is given by
\begin{equation*}
\begin{aligned}
\mathcal C^+
:= &\{\lambda\in \frak a_{\mathbb C}^* \,|\, \mathrm{Re} \langle \lambda, \gamma^\vee\rangle>0, ~\forall \gamma>0\}\\
=&\{z_1(2\alpha+\beta)+z_2(\alpha+\beta)\,|\,\mathrm{Re}z_1>\mathrm{Re}z_2>0\}.
\end{aligned}
\end{equation*}

For a positive root $\gamma$, denote by $w_\gamma$ the reflection defined by $\gamma$,
i.e., the reflection on the space $\frak a_{\mathbb C}^*$ which reflects $\gamma$ to $-\gamma$.
And denote by $\sigma(\omega)$ the rotation through $\omega$ with center at the origin.
Then it is well known that the Weyl group of $G_2$ is given by
\begin{equation*}
W =
\Big\{
e, ~w_{\alpha}, ~w_{\beta}, ~w_{3\alpha+\beta}, ~w_{2\alpha+\beta}, ~w_{3\alpha+2\beta}, ~w_{\alpha+\beta},
~\sigma(\frac{\pi}{3}), ~\sigma(\frac{2\pi}{3}), ~\sigma({\pi}), ~\sigma(\frac{4\pi}{3}), ~\sigma(\frac{5\pi}{3})
\Big\}.
\end{equation*}
Moreover by a direct calculation, we have the following table on $w\lambda$ and
\newline $\{\gamma>0\,|\,w\gamma<0\}$.
\begin{equation*}
\begin{matrix}
&w\lambda;\lambda=(z_1,z_2)&|&\{\gamma>0\,|\,w\gamma<0\}\\
e&(z_1,z_2)&|&-\\
w_\alpha&(z_2,z_1)&|&\alpha\\
w_\beta&(z_1+z_2,-z_2)&|&\beta\\
w_{3\alpha+\beta}&(-z_1,z_1+z_2)&|&\alpha,3\alpha+\beta,2\alpha+\beta\\
w_{2\alpha+\beta}&(-z_1-z_2,z_2)&|&\alpha,3\alpha+\beta,2\alpha+\beta,3\alpha+2\beta,\alpha+\beta\\
w_{3\alpha+2\beta}&(--z_2,-z_1)&|&3\alpha+\beta,2\alpha+\beta,3\alpha+2\beta,\alpha+\beta,\beta\\
w_{\alpha+\beta}&(z_1,--z_1-z_2)&|&3\alpha+2\beta,\alpha+\beta,\beta\\
\sigma(\frac{\pi}{3})&(-z_2,z_1+z_2)&|&\alpha+\beta,\beta\\
\sigma(\frac{2\pi}{3})&(-z_1-z_2,z_1)&|&2\alpha+\beta,3\alpha+2\beta,\alpha+\beta,\beta\\
\sigma(\pi)&(-z_1,-z_2)&|&\alpha, 3\alpha+\beta, 2\alpha+\beta,3\alpha+2\beta,\alpha+\beta,\beta\\
\sigma(\frac{4\pi}{3})&(z_2,-z_1-z_2)&|&\alpha, 3\alpha+\beta, 2\alpha+\beta,3\alpha+2\beta\\
\sigma(\frac{5\pi}{3})&(z_1+z_2,-z_1)&|&\alpha, 3\alpha+\beta
\end{matrix}
\end{equation*}
Also, by definition, we see that
\begin{equation*}
\begin{aligned}
&\langle\lambda,\alpha^\vee\rangle=z_1-z_2,\quad\qquad
\langle\lambda,\beta^\vee\rangle=z_2,\quad\qquad
\langle\lambda,(3\alpha+\beta)^\vee\rangle=z_1,\\
&\langle\lambda,(2\alpha+\beta)^\vee\rangle=2z_1+z_2,\
\langle\lambda,(3\alpha+2\beta)^\vee\rangle=z_1+z_2,\
\langle\lambda,(\alpha+\beta)^\vee\rangle=z_1+2z_2.
\end{aligned}
\end{equation*}
for $\lambda=(z_1,z_2)$,
since
\begin{equation*}
\begin{aligned}
&\lambda(t(x,x^{-1}))=x^{z_1}x^{-z_2}=x^{z_1-z_2},\quad
\lambda(t(1,x))=1^{z_1}x^{z_2}=x^{z_2},\\
&\lambda(t(x,1))=x^{z_1}1^{z_2}=x^{z_1},\quad
\lambda(t(x^2,x))=x^{2z_1}x^{z_2}=x^{2z_1+z_2},\\
&\lambda(t(x,x))=x^{z_1}x^{z_2}=x^{z_1+z_2},\quad
\lambda(t(x,x^{2}))=x^{z_1}x^{2z_2}=x^{z_1+2z_2}.
\end{aligned}
\end{equation*}
Hence, by tedious elementary calculations, which we decide to omit, we have the follows;
\vskip 0.30cm
\noindent
a) for $\displaystyle{ \langle w\lambda, \alpha^\vee \rangle-1}$
   and $\displaystyle{ \langle w\lambda, \beta^\vee \rangle-1}$,
\begin{equation*}
\begin{matrix}
&|&w\lambda;\lambda=(z_1,z_2)&|&\langle w\lambda,\alpha^\vee\rangle-1&|&\langle w\lambda,\beta^\vee\rangle-1\\
e&|&(z_1,z_2)&|&z_1-z_2-1&|&z_2-1\\
w_\alpha&|&(z_2,z_1)&|&z_2-z_1-1&|&z_1-1\\
w_\beta&|&(z_1+z_2,-z_2)&|&z_1+2z_2-1&|&-z_2-1\\
w_{3\alpha+\beta}&|&(-z_1,z_1+z_2)&|&-2z_1-z_2-1&|&z_1+z_2-1\\
w_{2\alpha+\beta}&|&(-z_1-z_2,z_2)&|&-z_1-2z_2-1&|&z_2-1\\
w_{3\alpha+2\beta}&|&(--z_2,-z_1)&|&z_1-z_2-1&|&-z_1-1\\
w_{\alpha+\beta}&|&(z_1,--z_1-z_2)&|&2z_1+z_2-1&|&-z_1-z_2-1\\
\sigma(\frac{\pi}{3})&|&(-z_2,z_1+z_2)&|&-z_1-2z_2-1&|&z_1+z_2-1\\
\sigma(\frac{2\pi}{3})&|&(-z_1-z_2,z_1)&|&-2z_1-z_2-1&|&z_1-1\\
\sigma(\pi)&|&(-z_1,-z_2)&|&-z_1+z_2-1&|&-z_2-1\\
\sigma(\frac{4\pi}{3})&|&(z_2,-z_1-z_2)&|&z_1+2z_2-1&|&-z_1-z_2-1\\
\sigma(\frac{5\pi}{3})&|&(z_1+z_2,-z_1)&|&2z_1+z_2-1&|&-z_1-1
\end{matrix}
\end{equation*}
and

\noindent
b) for $\displaystyle{\prod_{\gamma>0,w\gamma<0}\frac{\xi(\langle\lambda,\gamma^\vee\rangle)}{\xi(\langle\lambda,\gamma^\vee\rangle+1)}}$,
\begin{equation*}
\begin{matrix}
&\prod_{\gamma>0,w\gamma<0}\frac{\xi(\langle\lambda,\gamma^\vee\rangle)}{\xi(\langle\lambda,\gamma^\vee\rangle+1)}\\
e&1\\
w_\alpha&\frac{\xi(z_1-z_2)}{\xi(z_1-z_2+1)}\\
w_\beta&\frac{\xi(z_2)}{\xi(z_2+1)}\\
w_{3\alpha+\beta}&\frac{\xi(z_1-z_2)}{\xi(z_1-z_2+1)}\frac{\xi(z_1)}{\xi(z_1+1)}\frac{\xi(2z_1+z_2)}{\xi(2z_1+z_2+1)}\\
w_{2\alpha+\beta}&\frac{\xi(z_1-z_2)}{\xi(z_1-z_2+1)}\frac{\xi(z_1)}{\xi(z_1+1)}\frac{\xi(2z_1+z_2)}{\xi(2z_1+z_2+1)}\frac{\xi(z_1+z_2)}{\xi(z_1+z_2+1)}\frac{\xi(z_1+2z_2)}{\xi(z_1+2z_2+1)}\\
w_{3\alpha+2\beta}&\frac{\xi(z_1)}{\xi(z_1+1)}\frac{\xi(2z_1+z_2)}{\xi(2z_1+z_2+1)}\frac{\xi(z_1+z_2)}{\xi(z_1+z_2+1)}\frac{\xi(z_1+2z_2)}{\xi(z_1+2z_2+1)}\frac{\xi(z_2)}{\xi(z_2+1)}\\
w_{\alpha+\beta}&\frac{\xi(z_1+z_2)}{\xi(z_1+z_2+1)}\frac{\xi(z_1+2z_2)}{\xi(z_1+2z_2+1)}\frac{\xi(z_2)}{\xi(z_2+1)}\\
\sigma(\frac{\pi}{3})&\frac{\xi(z_1+2z_2)}{\xi(z_1+2z_2+1)}\frac{\xi(z_2)}{\xi(z_2+1)}\\
\sigma(\frac{2\pi}{3})&\frac{\xi(2z_1+z_2)}{\xi(2z_1+z_2+1)}\frac{\xi(z_1+z_2)}{\xi(z_1+z_2+1)}\frac{\xi(z_1+2z_2)}{\xi(z_1+2z_2+1)}\frac{\xi(z_2)}{\xi(z_2+1)}\\
\sigma(\pi)&\frac{\xi(z_1-z_2)}{\xi(z_1-z_2+1)}\frac{\xi(z_1)}{\xi(z_1+1)}\frac{\xi(2z_1+z_2)}{\xi(2z_1+z_2+1)}\frac{\xi(z_1+z_2)}{\xi(z_1+z_2+1)}\frac{\xi(z_1+2z_2)}{\xi(z_1+2z_2+1)}\frac{\xi(z_2)}{\xi(z_2+1)}\\
\sigma(\frac{4\pi}{3})&\frac{\xi(z_1-z_2)}{\xi(z_1-z_2+1)}\frac{\xi(z_1)}{\xi(z_1+1)}\frac{\xi(2z_1+z_2)}{\xi(2z_1+z_2+1)}\frac{\xi(z_1+z_2)}{\xi(z_1+z_2+1)}\\
\sigma(\frac{5\pi}{3})&\frac{\xi(z_1-z_2)}{\xi(z_1-z_2+1)}\frac{\xi(z_1)}{\xi(z_1+1)}
\end{matrix}
\end{equation*}
Or put them in a better form, we have
\begin{equation*}
\begin{matrix}
&\frac{1}{ \langle w\lambda,\alpha^\vee\rangle-1}\frac{1}{ \langle w\lambda,\beta^\vee\rangle-1}\cdot\prod_{\gamma>0,w\gamma<0}\frac{\xi(\langle\lambda,\gamma^\vee\rangle)}{\xi(\langle\lambda,\gamma^\vee\rangle+1)}\\
e&\frac{1}{z_1-z_2-1}\frac{1}{z_2-1}\\
w_\alpha&\frac{1}{z_2-z_1-1}\frac{1}{z_1-1}\cdot\frac{\xi(z_1-z_2)}{\xi(z_1-z_2+1)}\\
w_\beta&\frac{1}{z_1+2z_2-1}\frac{1}{-z_2-1}\cdot\frac{\xi(z_2)}{\xi(z_2+1)}\\
w_{3\alpha+\beta}&\frac{1}{-2z_1-z_2-1}\frac{1}{z_1+z_2-1}\cdot\frac{\xi(z_1-z_2)}{\xi(z_1-z_2+1)}\frac{\xi(z_1)}{\xi(z_1+1)}\frac{\xi(2z_1+z_2)}{\xi(2z_1+z_2+1)}\\
w_{2\alpha+\beta}&\frac{1}{-z_1-2z_2-1}\frac{1}{z_2-1}\cdot\frac{\xi(z_1-z_2)}{\xi(z_1-z_2+1)}\frac{\xi(z_1)}{\xi(z_1+1)}\frac{\xi(2z_1+z_2)}{\xi(2z_1+z_2+1)}\frac{\xi(z_1+z_2)}{\xi(z_1+z_2+1)}\frac{\xi(z_1+2z_2)}{\xi(z_1+2z_2+1)}\\
w_{3\alpha+2\beta}&\frac{1}{z_1-z_2-1}\frac{1}{-z_1-1}\cdot\frac{\xi(z_1)}{\xi(z_1+1)}\frac{\xi(2z_1+z_2)}{\xi(2z_1+z_2+1)}\frac{\xi(z_1+z_2)}{\xi(z_1+z_2+1)}\frac{\xi(z_1+2z_2)}{\xi(z_1+2z_2+1)}\frac{\xi(z_2)}{\xi(z_2+1)}\\
w_{\alpha+\beta}&\frac{1}{2z_1+z_2-1}\frac{1}{-z_1-z_2-1}\cdot\frac{\xi(z_1+z_2)}{\xi(z_1+z_2+1)}\frac{\xi(z_1+2z_2)}{\xi(z_1+2z_2+1)}\frac{\xi(z_2)}{\xi(z_2+1)}\\
\sigma(\frac{\pi}{3})&\frac{1}{-z_1-2z_2-1}\frac{1}{z_1+z_2-1}\cdot\frac{\xi(z_1+2z_2)}{\xi(z_1+2z_2+1)}\frac{\xi(z_2)}{\xi(z_2+1)}\\
\sigma(\frac{2\pi}{3})&\frac{1}{-2z_1-z_2-1}\frac{1}{z_1-1}\cdot\frac{\xi(2z_1+z_2)}{\xi(2z_1+z_2+1)}\frac{\xi(z_1+z_2)}{\xi(z_1+z_2+1)}\frac{\xi(z_1+2z_2)}{\xi(z_1+2z_2+1)}\frac{\xi(z_2)}{\xi(z_2+1)}\\
\sigma(\pi)&\frac{1}{-z_1+z_2-1}\frac{1}{-z_2-1}\cdot\frac{\xi(z_1-z_2)}{\xi(z_1-z_2+1)}\frac{\xi(z_1)}{\xi(z_1+1)}\frac{\xi(2z_1+z_2)}{\xi(2z_1+z_2+1)}\frac{\xi(z_1+z_2)}{\xi(z_1+z_2+1)}\frac{\xi(z_1+2z_2)}{\xi(z_1+2z_2+1)}\frac{\xi(z_2)}{\xi(z_2+1)}\\
\sigma(\frac{4\pi}{3})&\frac{1}{z_1+2z_2-1}\frac{1}{-z_1-z_2-1}\cdot\frac{\xi(z_1-z_2)}{\xi(z_1-z_2+1)}\frac{\xi(z_1)}{\xi(z_1+1)}\frac{\xi(2z_1+z_2)}{\xi(2z_1+z_2+1)}\frac{\xi(z_1+z_2)}{\xi(z_1+z_2+1)}\\
\sigma(\frac{5\pi}{3})&\frac{1}{2z_1+z_2-1}\frac{1}{-z_1-1}\cdot\frac{\xi(z_1-z_2)}{\xi(z_1-z_2+1)}\frac{\xi(z_1)}{\xi(z_1+1)}
\end{matrix}
\end{equation*}
By taking summation for all terms appeared, we then obtain the period $\omega_{\mathbb Q}^{G_2}(z_1,z_2)$ for $G_2$ over $\mathbb Q$.

\subsection{Zetas for $G_2$ over $\mathbb Q$}
Motivated by our study of high rank zeta functions in \cite{W1,W3} and new type of zetas for $SL(n)$ and $Sp(2n)$ in \cite{W4,W5}, as described in the introduction,
we can obtain two zeta functions for $G_2$ over $\mathbb Q$ from the period $\omega_{\mathbb Q}^{G_2}(z_1,z_2)$,
by taking residues along singular hyperplanes corresponding to (two) maximal parabolic subgroups.
\bigskip

\noindent
{\bf a) The zeta for $G_2/P_{\mathrm{long}}$}.

Recall that $P_{\mathrm{long}}$ corresponds to $\{\alpha\}=\Delta_0\backslash\{\beta\}$. Consequently, from the period $\omega_{\mathbb Q}^{G_2}(z_1,z_2)$ of $G_2$ over $\mathbb Q$, in order to introduce a zeta function $\xi_{\mathbb Q}^{G_2/P_{\mathrm{long}}}(s)$ for $G_2/P_{\mathrm{long}}$, we first take the residue along with the singular hyperplane $z_1-z_2=1$ of $\omega_{\mathbb Q}^{G_2}(z_1,z_2)$,
corresponding to $\langle \lambda+\rho_0,\alpha^\vee\rangle=0$, and set $z_2=s$ (then $z_1=1+s$ and $z_2-z_1=-1, ~2z_1+z_2=3s+2, ~z_1+z_2=2s+1, ~z_1+2z_2=3s+1, ~z_1-1=s, z_2+1=s+1$).
In such a way, we get then the following (single variable) period $\omega_{\mathbb Q}^{G_2/P{\mathrm{long}}}(s)$ associated to $G_2/P_{\mathrm{long}}$ over $\mathbb Q$: 

\begin{equation*}
\begin{aligned}
\omega_{\mathbb Q}^{G_2/P_{\mathrm{long}}}(s):=&\frac{1}{s-1}+\frac{1}{-2}\frac{1}{s}\cdot\frac{1}{\xi(2)}+0+
\frac{1}{-3s-3}\frac{1}{2s}\cdot \frac{1}{\xi(2)}\frac{\xi(s+1)}{\xi(s+2)}\frac{\xi(3s+2)}{\xi(3s+3)}\\
&+
\frac{1}{-3s-2}\frac{1}{s-1}\cdot \frac{1}{\xi(2)}\frac{\xi(s+1)}{\xi(s+2)}\frac{\xi(3s+2)}{\xi(3s+3)}\frac{\xi(2s+1)}{\xi(2s+2)}\frac{\xi(3s+1)}{\xi(3s+2)}\\
&+\frac{1}{-s-2}\cdot \frac{\xi(s+1)}{\xi(s+2)}\frac{\xi(3s+2)}{\xi(3s+3)}\frac{\xi(2s+1)}{\xi(2s+2)}\frac{\xi(3s+1)}{\xi(3s+2)}
\frac{\xi(s)}{\xi(s+1)}+0+0+0\\
&+
\frac{1}{-2}\frac{1}{-s-1}\cdot \frac{1}{\xi(2)}\frac{\xi(s+1)}{\xi(s+2)}\frac{\xi(3s+2)}{\xi(3s+3)}\frac{\xi(2s+1)}{\xi(2s+2)}\frac{\xi(3s+1)}{\xi(3s+2)}\frac{\xi(s)}{\xi(s+1)}\\
&+
\frac{1}{3s}\frac{1}{-2s-2}\cdot \frac{1}{\xi(2)}\frac{\xi(s+1)}{\xi(s+2)}\frac{\xi(3s+2)}{\xi(3s+3)}\frac{\xi(2s+1)}{\xi(2s+2)}\\
&+\frac{1}{3s+1}\frac{1}{-s-2}\cdot \frac{1}{\xi(2)}
\frac{\xi(s+1)}{\xi(s+2)}.
\end{aligned}
\end{equation*}
Multiplying with $\xi(2)\cdot\xi(s+2)\xi(2s+2)\xi(3s+3)$,
we then get
\begin{equation*}
\begin{aligned}
\xi_{\mathbb Q,o}^{G_2/P_{\mathrm{long}}}(s)=&
\frac{1}{s-1}\xi(2)\cdot\xi(s+2)\xi(2s+2)\xi(3s+3)\\
&-\frac{1}{s+2}\xi(2)\cdot\xi(s)\xi(2s+1)\xi(3s+1)\\
&-\frac{1}{2s}\cdot\xi(s+2)\xi(2s+2)\xi(3s+3)\\
&+\frac{1}{2(s+1)}\cdot\xi(s)\xi(2s+1)\xi(3s+1)\\
&-\frac{1}{3s+3}\frac{1}{2s}\cdot\xi(s+1)\xi(2s+2)\xi(3s+2)\\
&-\frac{1}{3s}\frac{1}{2s+2}\cdot\xi(s+1)\xi(2s+1)\xi(3s+2)\\
&-\frac{1}{3s+2}\frac{1}{s-1}\cdot\xi(s+1)\xi(2s+1)\xi(3s+1)\\
&-\frac{1}{3s+1}\frac{1}{s+2}\cdot\xi(s+1)\xi(2s+2)\xi(3s+3).
\end{aligned}
\end{equation*}
One checks easily the functional equation $\xi_{\mathbb Q,o}^{G_2/P_{\mathrm{long}}}(-1-s)=\xi_{\mathbb Q,o}^{G_2/P_{\mathrm{long}}}(s)$.
Define the first zeta function $\xi_{\mathbb Q}^{G_2/P_{\mathrm{long}}}(s)$
by normalizing $\xi_{\mathbb Q,o}^{G_2/P_{\mathrm{long}}}(s)$ with a shift
\begin{equation*}
\xi_{\mathbb Q}^{G_2/P_{\mathrm{long}}}(s):=\xi_{\mathbb Q,o}^{G_2/P_{\mathrm{long}}}(s-1).
\end{equation*}
Then we have the following
%
%%%%%%%%%%%%%%%%%%%%%%%%%%%%%%%%%%%%%%%%%%%%%%%%%%%%%%%%%%%%%%%%%%%%%%%%%%%%%%%%%%%%
%%
%% definition and proposition 1
%%
%%%%%%%%%%%%%%%%%%%%%%%%%%%%%%%%%%%%%%%%%%%%%%%%%%%%%%%%%%%%%%%%%%%%%%%%%%%%%%%%%%%%
%
\begin{defprop}
The zeta function $\xi_{\mathbb Q}^{G_2/P_{\mathrm{long}}}(s)$ for $(G_2, P_{\mathrm{long}})$ over $\mathbb Q$ given by
\begin{equation*}
\begin{aligned}
\xi^{G_2/P_{\mathrm{long}}}_{\mathbb Q}(s):=&\frac{1}{s-2}\xi(2)\cdot\xi(s+1)\xi(2s)\xi(3s)\\
&-\frac{1}{s+1}\xi(2)\cdot\xi(s-1)\xi(2s-1)\xi(3s-2)\\
&-\frac{1}{2s-2}\cdot\xi(s+1)\xi(2s)\xi(3s)\\
&+\frac{1}{2s}\cdot\xi(s-1)\xi(2s-1)\xi(3s-2)\\
&-\frac{1}{(3s)(2s-2)}\cdot\xi(s)\xi(2s)\xi(3s-1)\\
&-\frac{1}{(3s-1)(s-2)}\cdot\xi(s)\xi(2s-1)\xi(3s-2)\\
&-\frac{1}{(3s-3)(2s)}\cdot\xi(s)\xi(2s-1)\xi(3s-1)\\
&-\frac{1}{(3s-2)(s+1)}\cdot\xi(s)\xi(2s)\xi(3s),
\end{aligned}
\end{equation*}
satisfies the standard functional equation
\begin{equation*}
\xi_{\mathbb Q}^{G_2/P_{\mathrm{long}}}(1-s)=\xi_{\mathbb Q}^{G_2/P_{\mathrm{long}}}(s).
\end{equation*}
All poles of $\xi_{\Bbb Q}^{G_2/P_{\mathrm{long}}}(s)$ are two simple poles $s=-1,2$ and two double poles $s=0,1$.
\end{defprop}

\noindent
{\bf b) The zeta for $G_2/P_{\mathrm{short}}$}.
In parallel, recall that $P_{\mathrm{short}}$ corresponds to $\{\beta\}=\Delta_0\backslash\{\alpha\}$. Consequently, from the period $\omega_{\mathbb Q}^{G_2}(z_1,z_2)$ of $G_2$ over $\mathbb Q$, in order to introduce a zeta function $\xi_{\mathbb Q}^{G_2/P_{\mathrm{short}}}(s)$ for $G_2/P_{\mathrm{short}}$,  take the residue along  $z_2=1$,
corresponding to $\langle \lambda+\rho_0,\beta^\vee\rangle=0$,
and set $z_1=s$. Then we get accordingly for the period
$\omega_{\mathbb Q}^{G_2/P_{\mathrm{short}}}$ the following contributions:
\begin{equation*}
\begin{aligned}
\xi_{\mathbb Q,o}^{G_2/P_{\mathrm{short}}}(s)&:=\mathrm{Res}_{\langle\lambda+\rho_0,\beta^\vee\rangle=0}\omega_{\mathbb Q}^{G_2/P_2}(z_1,z_2):=\frac{1}{s-2}+0+\frac{1}{s+1}\frac{1}{-2}\cdot\frac{1}{\xi(2)}+0\\
&+\frac{1}{-s-3}\cdot\frac{\xi(s-1)}{\xi(s)}\frac{\xi(s)}{\xi(s+1)}\frac{\xi(2s+1)}{\xi(2s+2)}\frac{\xi(s+1)}{\xi(s+2)}\frac{\xi(s+2)}{\xi(s+3)}\\
&+
\frac{1}{s-2}\frac{1}{-s-1}\cdot \frac{\xi(s)}{\xi(s+1)}\frac{\xi(2s+1)}{\xi(2s+2)}\frac{\xi(s+1)}{\xi(s+2)}\frac{\xi(s+2)}{\xi(s+3)}\frac{1}{\xi(2)}\\
&+
\frac{1}{2s}\frac{1}{-s-2}\cdot\frac{\xi(s+1)}{\xi(s+2)}\frac{\xi(s+2)}{\xi(s+3)}\frac{1}{\xi(2)}\\
&+
\frac{1}{-s-3}\frac{1}{s}\cdot\frac{\xi(s+2)}{\xi(s+3)}\frac{1}{\xi(2)}\\
&+
\frac{1}{-s}\frac{1}{-2}\cdot \frac{\xi(s-1)}{\xi(s)}
\frac{\xi(s)}{\xi(s+1)}\frac{\xi(2s+1)}{\xi(2s+2)}\frac{\xi(s+1)}{\xi(s+2)}\frac{\xi(s+2)}{\xi(s+3)}\frac{1}{\xi(2)}
+0+0.
\end{aligned}
\end{equation*}
Multiplying with $\xi(2)\cdot\xi(s+3)\xi(2s+2)$,
and shift from $s$ to $s-1$ we then arrive at the second zeta function
$\xi^{G_2/P_{\mathrm{short}}}_{\mathbb Q}(s)$ for $(G_2, P_{\mathrm{short}})$ over $\mathbb Q$.
%
%%%%%%%%%%%%%%%%%%%%%%%%%%%%%%%%%%%%%%%%%%%%%%%%%%%%%%%%%%%%%%%%%%%%%%%%%%%%%%%%%%%%
%%
%% definition and proposition 2
%%
%%%%%%%%%%%%%%%%%%%%%%%%%%%%%%%%%%%%%%%%%%%%%%%%%%%%%%%%%%%%%%%%%%%%%%%%%%%%%%%%%%%%
%
\begin{defprop}
The zeta function $\xi^{G_2/P_{\mathrm{short}}}_{\mathbb Q}(s)$ for $(G_2, P_{\mathrm{short}})$ over $\mathbb Q$ given by
\begin{equation*}
\begin{aligned}
\xi^{G_2/P_{\mathrm{short}}}_{\mathbb Q}(s)=&\frac{1}{s-3}\xi(2)\cdot\xi(s+2)\xi(2s)\\
&-\frac{1}{s+2}\xi(2)\cdot\xi(s-2)\xi(2s-1)\\
&+\frac{1}{2s-2}\cdot\xi(s-2)\xi(2s-1)\\
&-\frac{1}{2s}\cdot\xi(s+2)\xi(2s)\\
&-\frac{1}{s(s-3)}\cdot\xi(s-1)\xi(2s-1)\\
&-\frac{1}{(s-1)(s+2)}\cdot\xi(s+1)\xi(2s)\\
&-\frac{1}{(2s-2)(s+1)}\cdot\xi(s)\xi(2s)\\
&-\frac{1}{(2s)(s-2)}\cdot\xi(s)\xi(2s-1),
\end{aligned}
\end{equation*}
satisfies the standard functional equation
\begin{equation*}
\xi^{G_2/P_{\mathrm{short}}}_{\mathbb Q}(1-s)=\xi^{G_2/P_{\mathrm{short}}}_{\mathbb Q}(s).
\end{equation*}
All poles of $\xi_{\Bbb Q}^{G_2/P_{\mathrm{short}}}(s)$ are four simple poles $s=-2,0,1,3$.
\end{defprop}

We expect that $\xi_{\mathbb Q}^{G_2/P_{\mathrm{long}}}(s)$ and $\xi_{\mathbb Q}^{G_2/P_{\mathrm{short}}}(s)$ satisfy the RH. For this, we have the following:
%
%%%%%%%%%%%%%%%%%%%%%%%%%%%%%%%%%%%%%%%%%%%%%%%%%%%%%%%%%%%%%%%%%%%%%%%%%%%%%%%%%%%%
%%
%% Theorem 5
%%
%%%%%%%%%%%%%%%%%%%%%%%%%%%%%%%%%%%%%%%%%%%%%%%%%%%%%%%%%%%%%%%%%%%%%%%%%%%%%%%%%%%%
%
\begin{thm}[Riemann Hypothesis$_{\mathbb Q}^{G_2/P}$]
\begin{equation*}
All\ zeros\ of\ \xi_{\mathbb Q}^{G_2/P_{\mathrm{long}}}(s)\ and\ \xi_{\mathbb Q}^{G_2/P_{\mathrm{short}}}(s)\ lie\ on\ the\ central\ line\ \mathrm{Re}(s)=1/2.
\end{equation*}
\end{thm}
\bigskip
\noindent
{\it Remark.} Zetas $\xi_{\mathbb Q}^{G_2/P}(s)$ are special cases of a more general construction;
In \cite{W4,W5}, we are able to define zeta functions $\xi_{\mathbb Q}^{G/P}(s)$
associated to classical semi-simple groups $G$ and their maximal parabolic subgroups $P$.
In particular, the conjectural standard functional equation
and the RH have been checked for $G=SL(2), ~SL(3), ~Sp(4)$ (\cite{W4,W5} for the FE, \cite{LS},  \cite{S,S2} for the RH).
Also, numerical calculations made by MS
give supportive evidences for the RH when $G=SL(4)$ or $SL(5)$.

%
%%%%%%%%%%%%%%%%%%%%%%%%%%%%%%%%%%%%%%%%%%%%%%%%%%%%%%%%%%%%%%%%%%%%%%%%%%%%%%%%%%%%%%%%%%%%%%%%%%%%%%%%%%%%%%%%%%%%%%%%%%
%%
%% \section 4
%%
%%%%%%%%%%%%%%%%%%%%%%%%%%%%%%%%%%%%%%%%%%%%%%%%%%%%%%%%%%%%%%%%%%%%%%%%%%%%%%%%%%%%%%%%%%%%%%%%%%%%%%%%%%%%%%%%%%%%%%%%%%
%
\section{Proof of the RH for $G_2$. Preliminaries}

To prove the RH for $G_2$, we prepare several auxiliary entire functions.
First, we define
\begin{equation*}
Z_1(s) := 12 s^3 (s-1)^3 \cdot (s+1)(3s-1)(2s-1)(3s-2)(s-2) \cdot \xi_{\Bbb Q}^{G_2/P_1}(s)
\end{equation*}
and
\begin{equation*}
Z_2(s) := 4 s^2(s-1)^2 \cdot (s+2)(s+1)(2s-1)(s-2)(s-3) \cdot \xi_{\Bbb Q}^{G_2/P_2}(s).
\end{equation*}
(Here, we use the notation $P_{\mathrm{long}}=P_{\beta}=P_1$ and $P_{\mathrm{short}}=P_\alpha=P_2$.)
Then $Z_1(s)$ and $Z_2(s)$ are entire functions by the results of section 3.
We have
\begin{equation*}
\aligned
Z_1(s) & = (s-1)\chi(2s) \Bigl[
(s-1)(3s-2)(As-A+1)\chi(s+1)\chi(3s) \\
& \quad - (s+1)(s-2)\chi(s)\chi(3s-1) - 2(s-1)(s-2)\chi(s)\chi(3s) \Bigr] \\
& \quad - s \, \chi(2s-1) \Bigl[ s(3s-1)(As-1)\chi(s-1)\chi(3s-2) \\
& \quad + (s+1)(s-2)\chi(s)\chi(3s-1) + 2s(s+1)\chi(s)\chi(3s-2) \Bigr] ,
\endaligned
\end{equation*}
and
\begin{equation*}
\aligned
Z_2(s) & = (s-2)\chi(2s) \Bigl[ (As + 3)(s-1)^2\chi(s+2)  \\
& \quad - 2(s-1)(s-3)\chi(s+1) - (s+2)(s-3)\chi(s) \Bigr]  \\
& \quad - (s+1) \chi(2s-1) \Bigl[ (As - A - 3)s^2\chi(s-2)   \\
& \quad + 2s(s+2) \, \chi(s-1) + (s+2)(s-3)\chi(s) \Bigr],
\endaligned
\end{equation*}
where
\begin{equation*}
\aligned
A &= 2\xi(2)-1 = \pi/3 -1 > 0, \\
\chi(s) &= s(s-1)\xi(s) = s(s-1)\pi^{-s/2}\Gamma(s/2)\zeta(s).
\endaligned
\end{equation*}
We find that
\begin{enumerate}
\item[$\bullet$]  $Z_1(s)$ has real zeros at $s=0, \, 1/3, \, 2/3, \, 1$ and $s=1/2$,
because all poles of $\xi_{\Bbb Q}^{G_2/P_1}(s)$ are two simple poles $s=-1,2$
and two double poles  $s=0,1$.
\item[$\bullet$]  $Z_2(s)$ has real zeros at $s=-1, \, 0, \, 1, \, 2$ and $s=1/2$,
because all poles of $\xi_{\Bbb Q}^{G_2/P_2}(s)$ are four simple poles $s=-2,0,1,3$.
\end{enumerate}
Hence the following two theorems are equivalent to
the RH of $\xi_{\Bbb Q}^{G_2/P_1}(s)$ and $\xi_{\Bbb Q}^{G_2/P_2}(s)$ respectively.
%
%%%%%%%%%%%%%%%%%%%%%%%%%%%%%%%%%%%%%%%%%%%%%%%%%%%%%%%%%%%%%%%%%%%%%%%%%%%%%%%%%%%%
%%
%% Theorem 6
%%
%%%%%%%%%%%%%%%%%%%%%%%%%%%%%%%%%%%%%%%%%%%%%%%%%%%%%%%%%%%%%%%%%%%%%%%%%%%%%%%%%%%%
%
\begin{thm} \label{thm_401}
All zeros of $Z_1(s)$ lie on the line $Re(s)=1/2$
except for four simple zeros $s=0, \, 1/3, \, 2/3, \, 1$.
\end{thm}
%
%%%%%%%%%%%%%%%%%%%%%%%%%%%%%%%%%%%%%%%%%%%%%%%%%%%%%%%%%%%%%%%%%%%%%%%%%%%%%%%%%%%%
%%
%% Theorem 7
%%
%%%%%%%%%%%%%%%%%%%%%%%%%%%%%%%%%%%%%%%%%%%%%%%%%%%%%%%%%%%%%%%%%%%%%%%%%%%%%%%%%%%%
%
\begin{thm} \label{thm_402}
All zeros of $Z_2(s)$ lie on the line $Re(s)=1/2$
except for four simple zeros $s=-1, \, 0, \, 1, \, 2$.
\end{thm}
Now we define
\begin{equation*}
\aligned
\widetilde{f}_1(s) & =   (s-1)(3s-2)(As-A+1)\chi(s+1)\chi(3s) \\
       & \quad - (s+1)(s-2)\chi(s)\chi(3s-1) - 2(s-1)(s-2)\chi(s)\chi(3s),
\endaligned
\end{equation*}
\begin{equation*}
\aligned
\widetilde{f}_2(s) & =   (As + 3)(s-1)^2\chi(s+2) \\
& \quad - 2(s-1)(s-3)\chi(s+1) - (s+2)(s-3)\chi(s).
\endaligned
\end{equation*}
and
\begin{equation} \label{009}
f_1(s)=(s-1)\widetilde{f}_1(s), \quad f_2(s)=(s-2)\widetilde{f}_2(s).
\end{equation}
Then
\[
Z_1(s) = \chi(2s) f_1(s) - \chi(2s-1)f_1(1-s),
\]
\[
Z_2(s) = \chi(2s) f_2(s) - \chi(2s-1)f_2(1-s).
\]

The proofs of Theorem \ref{thm_401} and Theorem \ref{thm_402} are  are divided into two steps.
First, we prove that all zeros of $f_i(s)$ lie in a vertical strip $\sigma_0 < \Re(s) <0$
except for finitely many exceptional zeros (section 5).
Then we obtain a nice product formula of $f_i(s)$ by a variant of Lemma 3 in~\cite{S}
(Lemma \ref{lem_01} in section 5, it will be proved in section 7).
Second, by using the product formula of $f_i(s)$,
we prove that all zeros of $Z_i(s)$ lie on the line $\Re(s)=1/2$ except for two simple zeros (section 6).
In this process, we use the result of Lagarias~\cite{La99} concerning the explicit upper bound
for the difference of the imaginary parts of the zeros of the Riemann zeta function.  (See also \cite{S2}.)
\smallskip

Before the proof, we recall the following result

\begin{lem}[\cite{LS}]
Let $\xi(s)$ the completed Riemann zeta function and $\chi(s)=s(s-1)\xi(s)$.
Then we have
\begin{equation} \label{LaSu}
\Bigl| \frac{\chi(2s-1)}{\chi(2s)} \Bigr| < 1 \quad \text{for} \quad {\rm Re}(s)>\frac{1}{2}.
\end{equation}

\end{lem}

%
%%%%%%%%%%%%%%%%%%%%%%%%%%%%%%%%%%%%%%%%%%%%%%%%%%%%%%%%%%%%%%%%%%%%%%%%%%%%%%%%%%%%%%%%%%%%%%%%%%%%%%%%%%%%%%%%%%%%%%%%%%
%%
%% section 5
%%
%%%%%%%%%%%%%%%%%%%%%%%%%%%%%%%%%%%%%%%%%%%%%%%%%%%%%%%%%%%%%%%%%%%%%%%%%%%%%%%%%%%%%%%%%%%%%%%%%%%%%%%%%%%%%%%%%%%%%%%%%%
%
\section{Proof of the RH for $G_2$: first step}

The aim of this section is to prove the following proposition.
%
%%%%%%%%%%%%%%%%%%%%%%%%%%%%%%%%%%%%%%%%%%%%%%%%%%%%%%%%%%%%%%%%%%%%%%%%%%%%%%%%%%%%
%%
%% Proposition 3
%%
%%%%%%%%%%%%%%%%%%%%%%%%%%%%%%%%%%%%%%%%%%%%%%%%%%%%%%%%%%%%%%%%%%%%%%%%%%%%%%%%%%%%
%
\begin{prop} \label{prop_01}
Let $f_1(s)$ and $f_2(s)$ be functions defined in \eqref{009}.
Then $f_i(s)~(i=1,2)$ has the product formula
\begin{equation*}
f_i(s) = C_i \,s^{m_i} \, e^{B_i^\prime s}
\Bigl( 1-\frac{s}{\beta_{0,i}} \Bigl)
\Bigl( 1-\frac{s}{\rho_{0,i}} \Bigl)
\Bigl( 1-\frac{s}{\overline{\rho}_{0,i}} \Bigl) \cdot \Pi_i(s)
\quad (B_i^\prime \geq 0),
\end{equation*}
where $C_1=f_1^\prime(0)$, $C_2=f_2(0)$, $\beta_{0,1}=1$, $\beta_{0,2}=2$, $m_1=1$, $m_2=0$,
$\rho_{0,i}~(i=1,2)$ are a complex zero of $f_i(s)$ with $\Re(\rho_{0,i})>1/2$ and
\begin{equation*}
\Pi_i(s) = \prod_{{\beta_i < 1/2}\atop{0 \not= \beta_i \in {\Bbb R}}} \Bigl( 1 -\frac{s}{\beta_i} \Bigr)
\prod_{{\beta_i<1/2}\atop{\gamma_i >0}}
\left[
\Bigl( 1 -\frac{s}{\rho_i} \Bigr) \Bigl( 1 -\frac{s}{\overline{\rho_i}} \Bigr)
\right]
\quad (\rho_i=\beta_i+\sqrt{-1} \, \gamma_i).
\end{equation*}
Here $\beta_i$ are at most finitely many real zeros of $f_i(s)$
and $\rho_i=\beta_i+i\gamma_i$ are other complex zeros of $f_i(s)$.
The product $\Pi_i(s)$ converges absolutely on any compact subset of ${\C}$
if we taken the product with the bracket.
\end{prop}
To prove the proposition, we prepare the following lemma:
%
%%%%%%%%%%%%%%%%%%%%%%%%%%%%%%%%%%%%%%%%%%%%%%%%%%%%%%%%%%%%%%%%%%%%%%%%%%%%%%%%%%%%
%%
%% lemma 1
%%
%%%%%%%%%%%%%%%%%%%%%%%%%%%%%%%%%%%%%%%%%%%%%%%%%%%%%%%%%%%%%%%%%%%%%%%%%%%%%%%%%%%%
%
\begin{lem} \label{lem_01}
Let $F(s)$ be an entire function of genus zero or one.
Suppose that
\begin{enumerate}
\item[{\rm (i)}] $F(s)$ is real on the real axis,
\item[{\rm (ii)}] there exists $\sigma_0 >0$ such
that all zeros of $F(s)$ lie in the vertical strip
\begin{equation*}
\sigma_0 < \Re(s) < 1/2
\end{equation*}
except for finitely many zeros,
\item[{\rm (iii)}]
the zeros of $F(s)$ are finitely many in the right-half plance $\Re(s) \geq 1/2$,
\item[{\rm (iv)}] there exists $C>0$ such that
\begin{equation}\label{le02}
N(T) \leq C T \log T \quad \text{as} \quad T \to \infty,
\end{equation}
where $N(T)$ is the number of zeros of $F(s)$
satisfying $0 \leq \Im(\rho) <T$.
\item[{\rm (v)}] $F(1-\sigma)/F(\sigma) > 0$ for large $\sigma>0$ and
\begin{equation}\label{le03}
F(1-\sigma)/ F(\sigma) \to 0 \quad \text{as} \quad \sigma \to \infty.
\end{equation}
\end{enumerate}
Then $F(s)$ has the product formula
\begin{equation*}
F(s)= C s^m e^{B^\prime s}
\prod_{0 \not=\rho \in {\Bbb R}}\Bigl(1-\frac{s}{\rho}\Bigr)
\prod_{\Im(\rho)>0}
\left[ \Bigl(1-\frac{s}{\rho}\Bigr)\Bigl(1-\frac{s}{\bar{\rho}}\Bigr) \right]
\end{equation*}
with
\begin{equation*}
B^\prime \geq 0,
\end{equation*}
The product in the right-hand side converges absolutely on every compact set
if we taken the product with the bracket.
\end{lem}

The most important part of this lemma is nonnegativity of $B^\prime$.
We will prove Lemma \ref{lem_01} in section 7.

If $f_i(s)$ satisfies all conditions in Lemma \ref{lem_01},
then we obtain Proposition \ref{prop_01} by applying Lemma \ref{lem_01} to $f_i(s)$.
Condition (i) is trivial for $f_i(s)$.
Under condition (ii), (iv) is easily proved by a standard argument
by using well-known estimate $|\chi(s)| \leq \exp(C|s|\log|s|)$ and Jensen's formula (see \S4.1 of ~\cite{S}, for example).
On the other hand, we have
\begin{equation*}
\aligned
f_1(0) & = 0, \quad f_1(s)  = f_1^\prime(0) \, s + O(s^2), \quad f_1^\prime(0) \simeq -2.176 \not= 0,\\
f_2(0) & \simeq -6.283 \not=0.
\endaligned
\end{equation*}
Hence it remains to prove (ii), (iii) and (v) for $f_i(s)$.

\subsection{Proof of (v)}

\subsubsection{Case of $f_1(s)$}

First we see that $f_1(1-\sigma)/f_1(\sigma)$ is positive for sufficiently large $\sigma>0$.
Using the functional equation of $\chi(s)$, we have
\begin{equation*}
\aligned
f_1(1-\sigma) & =
\sigma^2 (3\sigma-1)(A\sigma-1)\chi(\sigma-1)\chi(3\sigma-2) \\
& \quad + \sigma(\sigma-2)(\sigma+1) \chi(\sigma)\chi(3\sigma-1) + 2 \sigma^2 (\sigma+1)\chi(\sigma)\chi(3\sigma-2),
\endaligned
\end{equation*}
\begin{equation*}
\aligned
f_1(\sigma) & =
(\sigma-1)^2(3\sigma-2)(A\sigma-A+1)\chi(\sigma+1)\chi(3\sigma) \\
& \quad - (\sigma+1)(\sigma-2)(\sigma-1)\chi(\sigma)\chi(3\sigma-1) - 2(\sigma-1)^2(\sigma-2)\chi(\sigma)\chi(3\sigma)  .
\endaligned
\end{equation*}
Clearly the numerator is positive for large $\sigma>0$.
The denominator is also positive for large $\sigma>0$,
since $A>0$ and
\begin{equation} \label{015}
|\chi(\sigma)/\chi(\sigma+1)|<1 \quad (\sigma>0), \quad
|\chi(3\sigma-1)/\chi(3\sigma)|<1 \quad (\sigma>1/3)
\end{equation}
by replacing $2s-1$ by $\sigma$ or $3\sigma-1$ in \eqref{LaSu}.
Now we prove \eqref{le03}. We have
\[
\aligned
\frac{f_1(1-\sigma)}{f_1(\sigma)}
& =
\frac{\sigma^2 (3\sigma-1)(A\sigma-1)}{(\sigma-1)^2(3\sigma-2)(A\sigma-A+1)} \cdot
\frac{\chi(\sigma-1)\chi(3\sigma-2)}{\chi(\sigma+1)\chi(3\sigma)} \cdot
\frac{1+g(\sigma)}{1-h(\sigma)} \\
& =
\bigl( 1+O(\sigma^{-1}) \bigr) \cdot
\frac{\chi(\sigma-1)\chi(3\sigma-2)}{\chi(\sigma+1)\chi(3\sigma)} \cdot
\frac{1+g(\sigma)}{1-h(\sigma)},
\endaligned
\]
where
\[
g(\sigma)
=
\frac{(\sigma-2)(\sigma+1)}{\sigma (3\sigma-1)(A\sigma-1)} \cdot
\frac{\chi(\sigma)\chi(3\sigma-1)}{\chi(\sigma-1)\chi(3\sigma-2)}
+
\frac{2 (\sigma+1)}{(3\sigma-1)(A\sigma-1)} \cdot
\frac{\chi(\sigma)}{\chi(\sigma-1)},
\]
and
\[
h(\sigma)
=
\frac{(\sigma+1)(\sigma-2)}{(\sigma-1)(3\sigma-2)(A\sigma-A+1)} \cdot
\frac{\chi(\sigma)\chi(3\sigma-1)}{\chi(\sigma+1)\chi(3\sigma)}
+
\frac{2(\sigma-2)}{(3\sigma-2)(A\sigma-A+1)} \cdot
\frac{\chi(\sigma)}{\chi(\sigma+1)}.
\]
We have
\[
\aligned
\frac{\chi(\sigma-1)\chi(3\sigma-2)}{\chi(\sigma+1)\chi(3\sigma)}
& = \bigl(1+O(\sigma^{-1})\bigr)
\frac{\xi(\sigma-1)\xi(3\sigma-2)}{\xi(\sigma+1)\xi(3\sigma)} \\
& = \bigl(1+O(\sigma^{-1})\bigr) \cdot \pi^{2} \cdot
\frac{\Gamma((\sigma-1)/2)\Gamma((3\sigma-2)/2)}
{\Gamma((\sigma+1)/2) \Gamma(3\sigma/2)}
\frac{\zeta(\sigma-1)\zeta(3\sigma-2)}
{\zeta(\sigma+1)\zeta(3\sigma)}
\\
& = \bigl(1+O(\sigma^{-1})\bigr) \cdot
\frac{\Gamma((\sigma-1)/2)}
{\Gamma((\sigma+1)/2) (3\sigma-2)} \cdot O(1)
\endaligned
\]
for large $\sigma>0$. Using the Stirling formula
\[
\Gamma(z)=\sqrt{\frac{2\pi}{z}} \Bigl(\, \frac{z}{e} \,\Bigr)^z \bigl( 1 + O_\varepsilon(|z|^{-1})  \bigr)
\quad (|z| \geq 1, ~|{\rm arg}\,z|<\pi-\varepsilon ),
\]
we obtain
\begin{equation} \label{016}
\frac{\chi(\sigma-1)\chi(3\sigma-2)}{\chi(\sigma+1)\chi(3\sigma)}
= O(\sigma^{-2}) \quad \text{as} \quad \sigma \to +\infty.
\end{equation}
On the other hand, by using the Stirling formula again, we have
\begin{equation} \label{017}
g(\sigma)=O(1)+O(\sigma^{-1/2})=O(1) \quad \text{as} \quad \sigma \to +\infty.
\end{equation}
For $h(\sigma)$, by using \eqref{015}, we have
\begin{equation} \label{018}
h(\sigma)=O(\sigma^{-1}) \quad \text{as} \quad \sigma \to +\infty.
\end{equation}
From \eqref{016}, \eqref{017} and \eqref{018}, we obtain
\begin{equation*}
\frac{f_1(1-\sigma)}{f_1(\sigma)} =O(\sigma^{-2}) \quad \text{as} \quad \sigma \to +\infty.
\end{equation*}
This shows condition (v) for $f_1(s)$. \hfill $\Box$

\subsubsection{Case of $f_2(s)$}
First we see that $f_2(1-\sigma)/f_2(\sigma)$ is positive for sufficiently large $\sigma>0$.
Using the functional equation of $\chi(s)$, we have
\[
\aligned
f_2(1-\sigma) & =
\sigma^2 (\sigma+1)(A\sigma -A -  3)  \chi(\sigma-2) \\ & \quad
+ 2\sigma (\sigma+1)(\sigma+2)\chi(\sigma-1) + (\sigma-3)(\sigma+1)(\sigma+2)\chi(\sigma),
\endaligned
\]
\[
\aligned
f_2(\sigma) & =
(\sigma-1)^2 (\sigma-2)(A\sigma + 3) \chi(\sigma+2) \\ & \quad
- 2(\sigma-1)(\sigma-2)(\sigma-3)\chi(\sigma+1) - (\sigma+2)(\sigma-2)(\sigma-3)\chi(\sigma).
\endaligned
\]
Clearly the numerator is positive for large $\sigma>0$.
The denominator is also positive for large $\sigma>0$,
since $A>0$ and
\begin{equation} \label{020}
|\chi(\sigma+1)/\chi(\sigma+2)|<1 \quad (\sigma>-1), \quad
|\chi(\sigma)/\chi(\sigma+2)|<1 \quad (\sigma>0)
\end{equation}
by~\eqref{LaSu} and $\chi(\sigma)/\chi(\sigma+2)=(\chi(\sigma+1)/\chi(\sigma+2))\cdot(\chi(\sigma)/\chi(\sigma+1))$.
We have
\[
\aligned
\frac{f_2(1-\sigma)}{f_2(\sigma)}
& =
\frac{\sigma^2 (\sigma+1)(A\sigma -A -  3)}{(\sigma-1)^2 (\sigma-2)(A\sigma + 3)} \cdot
\frac{\chi(\sigma-2)}{\chi(\sigma+2)} \cdot
\frac{1+g(\sigma)}{1-h(\sigma)} \\
& =
\bigl( 1+O(\sigma^{-1}) \bigr) \cdot
\frac{\chi(\sigma-2)}{\chi(\sigma+2)} \cdot
\frac{1+g(\sigma)}{1-h(\sigma)},
\endaligned
\]
where
\[
g(\sigma)
=
\frac{2(\sigma+2)}{\sigma (A\sigma -A -  3)} \cdot
\frac{\chi(\sigma-1)}{\chi(\sigma-2)}
+
\frac{(\sigma+2)(\sigma-3)}{\sigma^2 (A\sigma -A -  3)} \cdot
\frac{\chi(\sigma)}{\chi(\sigma-2)},
\]
and
\[
h(\sigma)
=
\frac{2(\sigma-3)}{(\sigma-1)(A\sigma + 3)} \cdot
\frac{\chi(\sigma+1)}{ \chi(\sigma+2)}
+
\frac{(\sigma+2)(\sigma-3)}{(\sigma-1)^2 (A\sigma + 3)} \cdot
\frac{\chi(\sigma)}{\chi(\sigma+2)}.
\]
Using the Stirling formula,  we obtain
\begin{equation} \label{021}
\frac{\chi(\sigma-2)}{\chi(\sigma+2)}
= O(\sigma^{-2}) \quad \text{as} \quad \sigma \to +\infty.
\end{equation}
and
\begin{equation} \label{022}
g(\sigma)=O(\sigma^{-1/2})+O(1)=O(1) \quad \text{as} \quad \sigma \to +\infty.
\end{equation}
Using \eqref{020}, we have
\begin{equation} \label{023}
h(\sigma)=O(\sigma^{-1}) \quad \text{as} \quad \sigma \to +\infty.
\end{equation}
From \eqref{021}, \eqref{022} and \eqref{023}, we obtain
\begin{equation*}
\frac{f_2(1-\sigma)}{f_2(\sigma)} =O(\sigma^{-2}) \quad \text{as} \quad \sigma \to +\infty.
\end{equation*}
This shows condition (v) for $f_2(s)$. \hfill $\Box$

\subsection{Proof of (ii) and (iii)}
%
%%%%%%%%%%%%%%%%%%%%%%%%%%%%%%%%%%%%%%%%%%%%%%%%%%%%%%%%%%%%%%%%%%%%%%%%%%%%%%%%%%%%
%%
%% lemma 2
%%
%%%%%%%%%%%%%%%%%%%%%%%%%%%%%%%%%%%%%%%%%%%%%%%%%%%%%%%%%%%%%%%%%%%%%%%%%%%%%%%%%%%%
%
\begin{lem} \label{lem_02}
The entire function $f_1(s)$ has no zero in certain left-half plane $\Re(s)<\sigma_1$.
\end{lem}

\noindent
{\bf Proof.}
Assume $\sigma=\Re(s)<0$.
We have
\[
f_1(s)  = -(s+1)(s-1)(s-2)\chi(s)\chi(3s-1)
\left[\,
1 +  R_1(s) - R_2(s)
\,\right],
\]
where
\[
\aligned
R_1(s) & = 2 \, \frac{s-1}{s+1}
\cdot \frac{\chi(3s)}{\chi(3s-1)} \\
R_2(s) & = \frac{(s-1)(3s-2)(As-A+1)}{(s+1)(s-2)}
\cdot \frac{\chi(s+1)\chi(3s)}{\chi(s)\chi(3s-1)}.
\endaligned
\]
Clearly the factor $(s+1)(s-1)(s-2)\chi(s)\chi(3s-1)$ has no zero in the left-half plane
$\Re(s)< -1$.
Using the functional equation, we have
\[
\aligned
R_1(s)
& = 2 \,
\frac{s-1}{s+1} \frac{\chi(1-3s)}{\chi(2-3s)}
 =
\frac{6s(s-1)}{(3s-2)(s+1)}
\frac{\xi(1-3s)}{\xi(2-3s)} \\
& =
\frac{6\sqrt{\pi}s(s-1)}{(3s-2)(s+1)}
\frac{\Gamma((1-3s)/2)}{  \, \Gamma((2-3s)/2)}
\frac{\zeta(1-3s)}{\zeta(2-3s)}.
\endaligned
\]
Therefore
\[
|R_1(s)|
\leq 2\sqrt{\pi}
\left|
\frac{s(s-1)}{(s-\frac{2}{3})(s+1)}
\right|
\left|
\frac{\Gamma((1-3s)/2)}{  \, \Gamma((2-3s)/2)}
\right|
\zeta(1-3\sigma)\zeta(2-3\sigma).
\]
If $\sigma=\Re(s)<0$, $|{\rm arg}((1-3s)/2)| < \pi/2$ and $|{\rm arg}(2-3s)/2| < \pi/2$.
Hence we can apply the Stirling formula for $\Re(s)<0$, and then
\[
\left|
\frac{\Gamma((1-3s)/2)}{  \, \Gamma((2-3s)/2)}
\right|
= \sqrt{\frac{2}{3}} \, |s|^{-1/2}  (1+O(|s|^{-1})) \quad (\Re(s)<0).
\]
On the other hand
\[
\zeta(1-3\sigma)\zeta(2-3\sigma) \to 1 \quad (\sigma \to -\infty).
\]
Therefore
\begin{equation} \label{025}
|R_1(s)|
\leq \sqrt{\frac{8\pi}{3}} \cdot |s|^{-1/2} \cdot (1+O(|s|^{-1})),
\end{equation}
if $\sigma=\Re(s)<0$, and $|s|$, $|\sigma|$ are both large.

On the other hand, using the functional equation, we have
\[
\aligned
R_2(s)
& = \frac{(s-1)(3s-2)(As-A+1)}{(s+1)(s-2)}
\cdot \frac{\chi(-s)\chi(1-3s)}{\chi(1-s)\chi(2-3s)} \\
& =
\frac{3s(As-A+1)}{(s-2)}
\cdot \frac{\xi(-s)\xi(1-3s)}{\xi(1-s)\xi(2-3s)} \\
& =
\frac{3 \pi s(As-A+1)}{(s-2)}
\cdot \frac{\Gamma(-s/2)\Gamma((1-3s)/2)}{\Gamma((1-s)/2)\Gamma((2-3s)/2)}
\cdot \frac{\zeta(-s)\zeta(1-3s)}{\zeta(1-s)\zeta(2-3s)}.
\endaligned
\]
Therefore
\[
\aligned
|R_2(s)|
& \leq 3 \pi A
\left|
\frac{s(s-1+A^{-1})}{(s-2)}
\right|
\left|
\frac{\Gamma(-s/2)}{\Gamma((1-s)/2)}
\right|
\left|
\frac{\Gamma((1-3s)/2)}{\Gamma((2-3s)/2)}
\right|
\\
& \qquad \times
\zeta(-\sigma)\zeta(1-\sigma)\zeta(1-3\sigma)\zeta(2-3\sigma).
\endaligned
\]
If $\sigma=\Re(s)<0$, each argument of $-s/2$, $(1-s)/2$, $(1-3s)/2$ and $(2-3s)/2$ is less than $\pi/2$.
Hence we can apply the Stirling formula for $\Re(s)<0$, and then
\[
\left|
\frac{\Gamma(-s/2)}{\Gamma((1-s)/2)}
\right|
= \sqrt{2} \, |s|^{-1/2}  (1+O(|s|^{-1})) \quad (\Re(s)<0).
\]
\[
\left|
\frac{\Gamma((1-3s)/2)}{  \, \Gamma((2-3s)/2)}
\right|
= \sqrt{\frac{2}{3}} \, |s|^{-1/2}  (1+O(|s|^{-1})) \quad (\Re(s)<0).
\]
We have
\[
\zeta(-\sigma)\zeta(1-\sigma)\zeta(1-3\sigma)\zeta(2-3\sigma) \to 1 \quad (\sigma \to -\infty).
\]
Therefore
\begin{equation} \label{026}
|R_2(s)|
\leq  2 \sqrt{3} \pi A \cdot (1+O(|s|^{-1})),
\end{equation}
if $\sigma=\Re(s)<0$, and $|s|$, $-\sigma$ are both large.
Here
\[
 2 \sqrt{3} \pi A = 0.51364\dots
\]
Hence \eqref{025} and \eqref{026} implies Lemma \ref{lem_02}.
\hfill $\Box$
\smallskip
%
%%%%%%%%%%%%%%%%%%%%%%%%%%%%%%%%%%%%%%%%%%%%%%%%%%%%%%%%%%%%%%%%%%%%%%%%%%%%%%%%%%%%
%%
%% lemma 3
%%
%%%%%%%%%%%%%%%%%%%%%%%%%%%%%%%%%%%%%%%%%%%%%%%%%%%%%%%%%%%%%%%%%%%%%%%%%%%%%%%%%%%%
%
\begin{lem} \label{lem_03}
The entire function $f_2(s)$ has no zero in certain left-half plane $\Re(s)<\sigma_2$.
\end{lem}
\noindent
{\bf Proof.}
Assume $\sigma=\Re(s)<0$.
We have
\[
f_2(s)  = - (s+2)(s-2)(s-3)\chi(s)
\left[\,
1 +  R_1(s) - R_2(s)
\,\right],
\]
where
\[
R_1(s) = 2 \, \frac{(s-1)}{(s+2)}
\cdot \frac{\chi(s+1)}{\chi(s)}, \quad
R_2(s) = \frac{(As + 3)(s-1)^2}{(s+2)(s-3)}
\cdot \frac{\chi(s+2)}{\chi(s)}.
\]
Clearly the factor $(s+2)(s-2)(s-3)\chi(s)$ has no zero in the left-half plane
$\Re(s)< -2$.
Using the functional equation, we have
\[
R_1(s)
= 2 \, \frac{(s-1)}{(s+2)}
\cdot \frac{\chi(-s)}{\chi(1-s)}
 =
2 \, \frac{s+1}{s+2}
\frac{\xi(-s)}{\xi(1-s)}
 =
2 \sqrt{\pi} \, \frac{s+1}{s+2}
\frac{\Gamma(-s/2)}{  \, \Gamma((1-s)/2)}
\frac{\zeta(-s)}{\zeta(1-s)}.
\]
Therefore
\[
|R_1(s)|
\leq 2\sqrt{\pi}
\left|
\frac{s+1}{s+2}
\right|
\left|
\frac{\Gamma(-s/2)}{  \, \Gamma((1-s)/2)}
\right|
\zeta(-\sigma)\zeta(1-\sigma).
\]
If $\sigma=\Re(s)<0$, $|{\rm arg}(-s/2)| < \pi/2$ and $|{\rm arg}(1-s)/2| < \pi/2$.
Hence we can apply the Stirling formula for $\Re(s)<0$, and then
\[
\left|
\frac{\Gamma(-s/2)}{\Gamma((1-s)/2)}
\right|
= \sqrt{2} \, |s|^{-1/2}  (1+O(|s|^{-1})) \quad (\Re(s)<0).
\]
On the other hand $\zeta(-\sigma)\zeta(1-\sigma) \to 1$ as $\sigma \to -\infty$.
Therefore
\begin{equation} \label{025_2}
|R_1(s)|
\leq \sqrt{8\pi} \cdot |s|^{-1/2} \cdot (1+O(|s|^{-1})),
\end{equation}
if $\sigma=\Re(s)<0$, and $|s|$, $|\sigma|$ are both large.
On the other hand, using the functional equation, we have
\[
\aligned
R_2(s)
& = \frac{(As + 3)(s-1)^2}{(s+2)(s-3)}
\cdot \frac{\chi(-1-s)}{\chi(1-s)}
 =
\frac{(As + 3)(s-1)(s+1)}{s(s-3)}
\cdot \frac{\xi(-1-s)}{\xi(1-s)} \\
& = \pi \,
\frac{(As + 3)(s-1)(s+1)}{s(s-3)}
\cdot \frac{\Gamma((-1-s)/2)}{\Gamma((1-s)/2)}
\cdot \frac{\zeta(-1-s)}{\zeta(1-s)}.
\endaligned
\]
Therefore
\[
|R_2(s)|
\leq \pi A
\left|
\frac{(s + 3A^{-1})(s-1)(s+1)}{s(s-3)}
\right|
\left|
\frac{\Gamma((-1-s)/2)}{\Gamma((1-s)/2)}
\right|
\zeta(-1-\sigma)\zeta(1-\sigma).
\]
If $\sigma=\Re(s)<0$, both arguments of $(-1-s)/2$, $(1-s)/2$ are less than $\pi/2$.
Hence we can apply the Stirling formula for $\Re(s)<0$, and then
\[
\left|
\frac{\Gamma((-1-s)/2)}{\Gamma((1-s)/2)}
\right|
= 2 \, |s|^{-1}  (1+O(|s|^{-1})) \quad (\Re(s)<0).
\]
We have $\zeta(-1-\sigma)\zeta(1-\sigma) \to 1$ as $\sigma \to -\infty$.
Therefore
\begin{equation} \label{026_2}
|R_2(s)|
\leq  2 \pi A \cdot (1+O(|s|^{-1})),
\end{equation}
if $\sigma=\Re(s)<0$, and $|s|$, $|\sigma|$ are both large.
Here
\[
 2 \pi A = 0. 29655 \dots
\]
Hence \eqref{025_2} and \eqref{026_2} implies Lemma \ref{lem_03}.
\hfill $\Box$
%
%%%%%%%%%%%%%%%%%%%%%%%%%%%%%%%%%%%%%%%%%%%%%%%%%%%%%%%%%%%%%%%%%%%%%%%%%%%%%%%%%%%%
%%
%% lemma 4
%%
%%%%%%%%%%%%%%%%%%%%%%%%%%%%%%%%%%%%%%%%%%%%%%%%%%%%%%%%%%%%%%%%%%%%%%%%%%%%%%%%%%%%
%
\begin{lem} \label{lem_04}
The entire function $f_1(s)$ has only finitely many zeros in the right-half plane $\Re(s) > 1/3$.
In particular, the number of zeros of $f_1(s)$ in $\Re(s) \geq 1/2$ is finite.
\end{lem}

\noindent
{\bf Proof.}
We have
\begin{equation} \label{027}
f_1(s)  = (s-1)^2(3s-2)(As-A+1)\chi(s+1)\chi(3s)
\left[\,
1 -  Q_1(s) -Q_2(s)
\,\right],
\end{equation}
where
\begin{equation*}
Q_1(s) = \frac{(s+1)(s-2)}{(s-1)(3s-2)(As-A+1)} \cdot
\frac{\chi(s)\chi(3s-1)}{\chi(s+1)\chi(3s)}
\end{equation*}
\begin{equation*}
Q_2(s) =
\frac{2(s-2)}{(3s-2)(As-A+1)} \cdot
\frac{\chi(s)}{\chi(s+1)}.
\end{equation*}
The factor $(s-1)^2(3s-2)(As-A+1)\chi(s+1)\chi(3s)$ has no zero in $\Re(s)>1/3$
except for $s=2/3$ and $s=1$.
Replacing $2s-1$ by $3s-1$ or $s$ in \eqref{LaSu}, we obtain
\begin{equation} \label{030}
\left| \frac{\chi(3s-1)}{\chi(3s)} \right|  < 1 \quad \left( \Re(s) > \frac{1}{3} \right), \quad
\left| \frac{\chi(s)}{\chi(s+1)} \right|  < 1 \quad (\Re(s) > 0).
\end{equation}
Let $D_1$ be the region
\[
D_1:=
\left\{
s \in {\mathbb C} ~\left|~
\Re(s) \geq \frac{1}{3}, ~\left| \frac{(s+1)(s-2)}{(s-1)(3s-2)(As-A+1)} \right|
+ \left| \frac{2(s-2)}{(3s-2)(As-A+1)} \right|\geq 1
\right.\right\}.
\]
Then $f_1(s) \not=0$ if $s \not\in D_1$ and $\Re(s) \geq 1/3$,
because of \eqref{027} and \eqref{030}.
The region $D_1$ is bounded, since
\[
\left| \frac{(s+1)(s-2)}{(s-1)(3s-2)(As-A+1)} \right|
+ \left| \frac{2(s-2)}{(3s-2)(As-A+1)} \right| < 1
\]
for large $|s|$. Hence the number of zeros of $f_1(s)$ in $\Re(s) \geq 1/3$ is finite.
\hfill $\Box$
\bigskip
%
%%%%%%%%%%%%%%%%%%%%%%%%%%%%%%%%%%%%%%%%%%%%%%%%%%%%%%%%%%%%%%%%%%%%%%%%%%%%%%%%%%%%
%%
%% lemma 5
%%
%%%%%%%%%%%%%%%%%%%%%%%%%%%%%%%%%%%%%%%%%%%%%%%%%%%%%%%%%%%%%%%%%%%%%%%%%%%%%%%%%%%%
%
\begin{lem} \label{lem_05}
The entire function $f_2(s)$ has only finitely many zeros in the right-half plane $\Re(s) > 0$.
In particular, the number of zeros of $f_2(s)$ in $\Re(s) \geq 1/2$ is finite.
\end{lem}

\noindent
{\bf Proof.}
We have
\begin{equation} \label{027_2}
f_2(s)  = (As + 3)(s-1)^2(s-2)\chi(s+2)
\left[\,
1 -  Q_1(s) -Q_2(s)
\,\right],
\end{equation}
where
\begin{equation*}
Q_1(s) = \frac{2 (s-3)}{(As + 3)(s-1)} \cdot
\frac{\chi(s+1)}{\chi(s+2)}
\end{equation*}
\begin{equation*}
Q_2(s) =
\frac{(s+2)(s-3)}{(As + 3)(s-1)^2} \cdot
\frac{\chi(s)}{\chi(s+2)}.
\end{equation*}
The factor $(As + 3)(s-1)^2(s-2)\chi(s+2)$ has no zero in $\Re(s)>0$ except for $s=1$ and $s=2$.
Replacing $2s-1$ by $s+1$ or $s$ in \eqref{LaSu}, we obtain
\begin{equation} \label{030_2}
\aligned
\left| \frac{\chi(s+1)}{\chi(s+2)} \right|  & < 1 \quad \left( \Re(s) > -1 \right), \\
\left| \frac{\chi(s)}{\chi(s+2)} \right| &
= \left| \frac{\chi(s+1)}{\chi(s+2)} \right|
\left| \frac{\chi(s)}{\chi(s+1)} \right|
 < 1 \quad (\Re(s) > 0).
\endaligned
\end{equation}
Let $D_2$ be the region
\[
D_2:=
\left\{
s \in {\mathbb C} ~\left|~
\Re(s) \geq 0, ~\left| \frac{2 (s-3)}{(As + 3)(s-1)} \right|
+ \left| \frac{(s+2)(s-3)}{(As + 3)(s-1)^2 } \right|\geq 1
\right.\right\}.
\]
Then $f_2(s) \not=0$ if $s \not\in D_2$ and $\Re(s) \geq 0$,
because of \eqref{027_2} and \eqref{030_2}.
Clearly the region $D_2$ is bounded, the number of zeros of $f_2(s)$ in $\Re(s) \geq 0$ is finite.
\hfill $\Box$
\bigskip

\subsection{Proof of Proposition \ref{prop_01}}

By the results in section 5,1 and section 5.2, we can apply Lemma \ref{lem_01} to $f_i(s)~(i=1,2)$.
Hence the proof of Proposition \ref{prop_01} is completed by the following lemmas.
%
%%%%%%%%%%%%%%%%%%%%%%%%%%%%%%%%%%%%%%%%%%%%%%%%%%%%%%%%%%%%%%%%%%%%%%%%%%%%%%%%%%%%
%%
%% lemma 6
%%
%%%%%%%%%%%%%%%%%%%%%%%%%%%%%%%%%%%%%%%%%%%%%%%%%%%%%%%%%%%%%%%%%%%%%%%%%%%%%%%%%%%%
%
\begin{lem} \label{lem_06}
The number of zeros of $f_1(s)$ in $\Re(s) \geq 1/2$ is just three,
One of them is the real zero $s=1$,
and another two zeros are non-real zeros and conjugate each other.
The values of complex zeros are about $s \simeq 0.90 \pm i \cdot 2.09$.
\end{lem}
%
%%%%%%%%%%%%%%%%%%%%%%%%%%%%%%%%%%%%%%%%%%%%%%%%%%%%%%%%%%%%%%%%%%%%%%%%%%%%%%%%%%%%
%%
%% lemma 7
%%
%%%%%%%%%%%%%%%%%%%%%%%%%%%%%%%%%%%%%%%%%%%%%%%%%%%%%%%%%%%%%%%%%%%%%%%%%%%%%%%%%%%%
%
\begin{lem} \label{lem_07}
The number of zeros of $f_2(s)$ in $\Re(s) \geq 1/2$ is just three.
One of them is the real zero $s=2$,
and another two zeros are non-real zeros and conjugate each other.
The values of complex zeros are $s \simeq 1.17 \pm i \cdot 3.43$.
\end{lem}
\noindent
{\bf Proof of Lemma \ref{lem_06} and Lemma \ref{lem_07}.}
The domain $D_i \cap \{\Re(s) \leq 1/2\}$
is contained in the rectangle $R=[1/2,5] \times [-10,10]$,
where $D_i$ is the region in the proof of Lemma \ref{lem_06} or Lemma \ref{lem_07}.
Because of the argument principle, the number of zeros of $f(s)$ in $R$ is given by
\[
\frac{1}{2\pi i} \int_{\partial R} \frac{f_i^\prime}{f_i}(s) ds.
\]
In particular, the value of this integral is an integer.
Therefore we can check that the value of this integral is just three
by a computational way (for example Mathematica, Maple, PARI$/$GP, etc.).
Hence we conclude that $f_i(s)$ has just three zeros in the rectangle $R$.
One of them are trivial real zero of $f_1(s)$ (resp. $f_2(s)$) at $s=1$ (resp. $s=2$).
By suitable computational way,
we find an approximate value of the above two complex zeros of $f_1(s)$ (resp. $f_2(s)$)
are $s \simeq 0.90 \pm i \cdot 2.09$ (resp. $s \simeq 1.17 \pm i \cdot 3.43$).
\hfill $\Box$
%
%%%%%%%%%%%%%%%%%%%%%%%%%%%%%%%%%%%%%%%%%%%%%%%%%%%%%%%%%%%%%%%%%%%%%%%%%%%%%%%%%%%%%%%%%%%%%%%%%%%%%%%%%%%%%%%%%%%%%%%%%%
%%
%% section 6
%%
%%%%%%%%%%%%%%%%%%%%%%%%%%%%%%%%%%%%%%%%%%%%%%%%%%%%%%%%%%%%%%%%%%%%%%%%%%%%%%%%%%%%%%%%%%%%%%%%%%%%%%%%%%%%%%%%%%%%%%%%%%
%
\section{Proof of the RH for $G_2$: second step}

\subsection{Proof of  Theorem \ref{thm_401} and  Theorem \ref{thm_402}}

We have the following three assertions for $Z_1(s)$.
%
%%%%%%%%%%%%%%%%%%%%%%%%%%%%%%%%%%%%%%%%%%%%%%%%%%%%%%%%%%%%%%%%%%%%%%%%%%%%%%%%%%%%
%%
%% proposition 4
%%
%%%%%%%%%%%%%%%%%%%%%%%%%%%%%%%%%%%%%%%%%%%%%%%%%%%%%%%%%%%%%%%%%%%%%%%%%%%%%%%%%%%%
%
\begin{prop} \label{prop_02}
$Z_1(s)$ has no zero in the right-half plane $\Re(s) \geq 20$.
\end{prop}
%
%%%%%%%%%%%%%%%%%%%%%%%%%%%%%%%%%%%%%%%%%%%%%%%%%%%%%%%%%%%%%%%%%%%%%%%%%%%%%%%%%%%%
%%
%% proposition 5
%%
%%%%%%%%%%%%%%%%%%%%%%%%%%%%%%%%%%%%%%%%%%%%%%%%%%%%%%%%%%%%%%%%%%%%%%%%%%%%%%%%%%%%
%
\begin{prop} \label{prop_03}
$Z_1(s)$ has no zero in the region $1/2 < \sigma < 20$, $|t| \geq 25$.
\end{prop}
%
%%%%%%%%%%%%%%%%%%%%%%%%%%%%%%%%%%%%%%%%%%%%%%%%%%%%%%%%%%%%%%%%%%%%%%%%%%%%%%%%%%%%
%%
%% proposition 6
%%
%%%%%%%%%%%%%%%%%%%%%%%%%%%%%%%%%%%%%%%%%%%%%%%%%%%%%%%%%%%%%%%%%%%%%%%%%%%%%%%%%%%%
%
\begin{prop} \label{prop_04}
$Z_1(s)$ has only one simple zero $s=2/3,\,1$ in the region $1/2 < \sigma < 20$, $|t| \leq 25$.
\end{prop}
Then, as a consequence of these results and the functional equation of $Z_1(s)$,
all zeros of $Z_1(s)$ lie on the line $\Re(s)=1/2$
except for simple zeros $s=0,\, 1/3, \, 1/2, \, 2/3, \, 1$.

While we have the following three assertions for $Z_2(s)$.
%
%%%%%%%%%%%%%%%%%%%%%%%%%%%%%%%%%%%%%%%%%%%%%%%%%%%%%%%%%%%%%%%%%%%%%%%%%%%%%%%%%%%%
%%
%% proposition 7
%%
%%%%%%%%%%%%%%%%%%%%%%%%%%%%%%%%%%%%%%%%%%%%%%%%%%%%%%%%%%%%%%%%%%%%%%%%%%%%%%%%%%%%
%
\begin{prop} \label{prop_05}
$Z_2(s)$ has no zero in the right-half plane $\Re(s) \geq 20$.
\end{prop}
%
%%%%%%%%%%%%%%%%%%%%%%%%%%%%%%%%%%%%%%%%%%%%%%%%%%%%%%%%%%%%%%%%%%%%%%%%%%%%%%%%%%%%
%%
%% proposition 8
%%
%%%%%%%%%%%%%%%%%%%%%%%%%%%%%%%%%%%%%%%%%%%%%%%%%%%%%%%%%%%%%%%%%%%%%%%%%%%%%%%%%%%%
%
\begin{prop} \label{prop_06}
$Z_2(s)$ has no zero in the region $1/2 < \sigma < 20$, $|t| \geq 36$.
\end{prop}
%
%%%%%%%%%%%%%%%%%%%%%%%%%%%%%%%%%%%%%%%%%%%%%%%%%%%%%%%%%%%%%%%%%%%%%%%%%%%%%%%%%%%%
%%
%% proposition 9
%%
%%%%%%%%%%%%%%%%%%%%%%%%%%%%%%%%%%%%%%%%%%%%%%%%%%%%%%%%%%%%%%%%%%%%%%%%%%%%%%%%%%%%
%
\begin{prop} \label{prop_07}
$Z_2(s)$ has only one simple zero $s=1, \, 2$ in the region $1/2 < \sigma < 20$, $|t| \leq 36$.
\end{prop}
Then, as a consequence of these results and the functional equation of $Z_2(s)$,
all zeros of $Z_2(s)$ lie on the line $\Re(s)=1/2$
except for simple zeros $s=-1,\, 0, \, 1/2, \, 1, \, 2$. \hfill $\Box$
\bigskip

Hence it remains to prove the above six propositions. We carry out the proof of them in below.
The hardest part is the proof of Proposition \ref{prop_03} and \ref{prop_06}.
To prove Proposition \ref{prop_03} and \ref{prop_06}, we use the results in the first step
and a result of Lagarias~\cite{La99}.

\subsection{Proof of Proposition \ref{prop_02}}
We have
\begin{equation} \label{3_term}
\aligned
Z_1(s)  & = (s-1)^2(3s-2)(As-A+1)\chi(s+1)\chi(3s)\chi(2s) \\
& \quad \times (1 - R_1(s) - R_2(s) - R_3(s) + R_4(s)+R_5(s)),
\endaligned
\end{equation}
where
\begin{equation*}
\aligned
R_1(s) & =  \frac{ (s+1)(s-2) }{ (s-1)(3s-2)(As-A+1) }
\frac{\chi(s)}{\chi(s+1)} \frac{\chi(3s-1)}{\chi(3s)}, \\
R_2(s) & =  \frac{ 2(s-1)(s-2) }{ (s-1)(3s-2)(As-A+1) }\frac{\chi(s)}{\chi(s+1)}, \\
R_3(s) & =  \frac{ s^2(3s-1)(As-1) }{ (s-1)^2(3s-2)(As-A+1)}
\frac{\chi(s-1)}{\chi(s+1)}\frac{\chi(3s-2)}{\chi(3s)} \frac{\chi(2s-1)}{\chi(2s)}, \\
R_4(s) & =  \frac{ s(s+1)(s-2) }{ (s-1)^2(3s-2)(As-A+1) }
\frac{\chi(s)}{\chi(s+1)}\frac{\chi(3s-1)}{\chi(3s)} \frac{\chi(2s-1)}{\chi(2s)}, \\
R_5(s) & =  \frac{ 2s^2(s+1) }{ (s-1)^2(3s-2)(As-A+1) }
\frac{\chi(s)}{\chi(s+1)}\frac{\chi(3s-2)}{\chi(3s)} \frac{\chi(2s-1)}{\chi(2s)}.
\endaligned
\end{equation*}
Replacing $2s-1$ by $s$ or $3s-1$ in \eqref{LaSu}, we have 
\[
\left| \frac{\chi(s)}{\chi(s+1)} \right| < 1 ~(\Re(s)>0), \quad
\left| \frac{\chi(3s-1)}{\chi(3s)} \right| < 1 ~\Bigl(\Re(s)>\frac{1}{3}\Bigr).
\]
Moreover, replacing $2s-1$ by $s-1$ or $3s-2$ in \eqref{LaSu}, we have 
\[
\aligned
\left| \frac{\chi(s-1)}{\chi(s+1)} \right|
& = \left| \frac{\chi(s)}{\chi(s+1)} \right| \left| \frac{\chi(s-1)}{\chi(s)} \right|
< 1 \quad (\Re(s)>1), \\
\left| \frac{\chi(3s-2)}{\chi(3s)} \right|
& = \left| \frac{\chi(3s-1)}{\chi(3s)} \right| \left| \frac{\chi(3s-2)}{\chi(3s-1)} \right|
< 1 \quad \Bigl(\Re(s)>\frac{2}{3}\Bigr).
\endaligned
\]
Hence $|R_i(s)| \leq C_i |s|^{-1}~(i=1,2,4,5)$ for $\Re(s)>1$.
Applying the Stirling formula to $R_3(s)$,
we obtain $|R_3(s)|=(|s|^{-5/2})$ for $\Re(s)>1$ as $|s|\to \infty$ in the right-half plane.
%\begin{equation*}
%\aligned
%|R_1(s)| & =O(|s|^{-2}), \quad |R_2(s)|=(|s|^{-3/2}), \quad |R_3(s)|=(|s|^{-5/2}), \\
%|R_4(s)| & =(|s|^{-3/2}), \quad |R_5(s)|=(|s|^{-2}),
%\endaligned
%\end{equation*}
Therefore $Z_1(s) \not=0$ for some right-half plane $\Re(s) \geq \sigma_3$.
Using the monotone decreasing property of $\zeta(\sigma)$ as $\sigma \to +\infty$
and the effective version of Stirling's formula (\cite{Ol})
\begin{equation*}
\Gamma(s) = \Bigl( \frac{2\pi}{s} \Bigr)^{\frac{1}{2}} \Bigl( \frac{s}{e} \Bigr)^s
\Bigl\{ 1+\Theta \Bigl( \frac{1}{8|s|}\Bigr) \Bigr\} \quad (\Re(s)>1),
\end{equation*}
where the notation $f=\Theta(g)$ means $|f| \leq g$,
% \begin{equation}
% |\Gamma(s)| = \sqrt{2\pi} |s|^{\sigma-\frac{1}{2}} e^{-\sigma-t\arctan(\frac{t}{\sigma})}
% \Bigl\{ 1+\Theta \Bigl( \frac{1}{8|s|}\Bigr) \Bigr\},
% \quad s=\sigma+it,~\sigma>0
% \end{equation}
we have
\begin{equation*}
|R_1(s)| \leq 0.1, \quad |R_2(s)| \leq 0.3, \quad |R_3(s)| \leq 0.05 \quad |R_4(s)| \leq 0.1, \quad |R_5(s)| \leq 0.1
\end{equation*}
for $\Re(s) \geq 20$ (in fact, these bounds already hold for $\Re(s) \geq 10$).
These estimates imply $Z_1(s)\not=0$ for $\Re(s) \geq 20$ by \eqref{3_term},
since $(s-1)^2(3s-2)(As-A+1)\chi(s+1)\chi(3s)\chi(2s)$ has no zero in the right-half plane $\Re(s) \geq 20$.
\hfill $\Box$

%
%%%%%%%%%%%%%%%%%%%%%%%%%%%%%%%%%%%%%%%%%%%%%%%%%%%%%%%%%%%%%%%%%%%%%%%%%%%%%%%%%%%%
%%
%% Proof of prop_05
%%
%%%%%%%%%%%%%%%%%%%%%%%%%%%%%%%%%%%%%%%%%%%%%%%%%%%%%%%%%%%%%%%%%%%%%%%%%%%%%%%%%%%%
%
\subsection{Proof of Proposition \ref{prop_05}}
We have
\begin{equation} \label{3_term_2}
\aligned
Z_2(s) & = (s-1)^2(s-2)(As + 3) \chi(s+2)\chi(2s) \\
&  \qquad \times (1 - R_1(s) - R_2(s) + R_3(s) - R_4(s) - R_5(s)),
\endaligned
\end{equation}
where
\begin{equation*}
\aligned
R_1(s) & =  \frac{ 2(s-3) }{ (s-1)(As + 3) } \frac{\chi(s+1)}{\chi(s+2)}, \\
R_2(s) & =  \frac{ (s+2)(s-3) }{ (s-1)^2(As + 3) } \frac{\chi(s)}{\chi(s+2)}, \\
R_3(s) & =  \frac{ s^2(s+1)(A s - 3 -A) }{ (s-1)^2(s-2)(As + 3) }
\frac{\chi(s-2)}{\chi(s+2)} \frac{\chi(2s-1)}{\chi(2s)}, \\
R_4(s) & =  \frac{ 2s(s+1)(s+2) }{ (s-1)^2(s-2)(As + 3) }
\frac{\chi(s-1)}{\chi(s+2)} \frac{\chi(2s-1)}{\chi(2s)}, \\
R_5(s) & =  \frac{ (s-3)(s+1)(s+2) }{ (s-1)^2(s-2)(As + 3) }
\frac{\chi(s)}{\chi(s+2)} \frac{\chi(2s-1)}{\chi(2s)}.
\endaligned
\end{equation*}
Replacing $2s-1$ by $s-a$ ($a=-1,0,1,2$) in \eqref{LaSu}, we have
\[
\aligned
\left| \frac{\chi(s+1)}{\chi(s+2)} \right| & < 1,  \quad (\Re(s)>-1), \\
\left| \frac{\chi(s)}{\chi(s+2)} \right|
& =  \left| \frac{\chi(s+1)}{\chi(s+2)} \right| \left| \frac{\chi(s)}{\chi(s+1)} \right|  < 1,  \quad (\Re(s)>0), \\
\left| \frac{\chi(s-1)}{\chi(s+2)} \right|
& =
\left| \frac{\chi(s+1)}{\chi(s+2)} \right|
\left| \frac{\chi(s)}{\chi(s+1)} \right|
\left| \frac{\chi(s-1)}{\chi(s)} \right|  < 1,  \quad (\Re(s)>1), \\
\left| \frac{\chi(s-2)}{\chi(s+2)} \right|
& =
\left| \frac{\chi(s+1)}{\chi(s+2)} \right|
\left| \frac{\chi(s)}{\chi(s+1)} \right|
\left| \frac{\chi(s-1)}{\chi(s)} \right|
\left| \frac{\chi(s-2)}{\chi(s-1)} \right| < 1,  \quad (\Re(s)>2). \\
\endaligned
\]
Hence $|R_i(s)| \leq C_i |s|^{-1}~(i=1,2,4,5)$ for $\Re(s)>2$.
Applying the Stirling formula to $R_3(s)$,
we obtain $|R_3(s)|=(|s|^{-5/2})$ for $\Re(s) \gg 0$.
Therefore $Z_2(s) \not=0$ for some right-half plane $\Re(s) \geq \sigma_4$.
Using the monotone decreasing property of $\zeta(\sigma)$ as $\sigma \to +\infty$
and the effective version of Stirling's formula,
we have
\begin{equation*}
|R_1(s)| \leq 0.3, \quad |R_2(s)| \leq 0.13, \quad |R_3(s)| \leq 0.15 \quad |R_4(s)| \leq 0.2, \quad |R_5(s)| \leq 0.1
\end{equation*}
for $\Re(s) \geq 20$.
These estimates imply $Z_2(s)\not=0$ for $\Re(s) \geq 20$ by \eqref{3_term_2},
since $(s-1)^2(s-2)(As + 3) \chi(s+2)\chi(2s)$ has no zero in the right-half plane $\Re(s) \geq 20$.
\hfill $\Box$

%
%%%%%%%%%%%%%%%%%%%%%%%%%%%%%%%%%%%%%%%%%%%%%%%%%%%%%%%%%%%%%%%%%%%%%%%%%%%%%%%%%%%%
%%
%% Proof of Lemma
%%
%%%%%%%%%%%%%%%%%%%%%%%%%%%%%%%%%%%%%%%%%%%%%%%%%%%%%%%%%%%%%%%%%%%%%%%%%%%%%%%%%%%%
%
\subsection{Proof of Proposition \ref{prop_03}}
Let $\rho_0=\beta_0+ i \gamma_0~(\gamma_0>0)$ be the complex zero of $f_1(s)$ in Lemma \ref{lem_06}.
By Proposition \ref{prop_01}, $f_1(s)$ has the factorization
\begin{equation*} \label{17}
f_1(s) = f_1^\prime (0) \, e^{B_1^\prime s} \,
s(1-s)
\Bigl( 1-\frac{s}{\rho_0} \Bigl)
\Bigl( 1-\frac{s}{\overline{\rho}_0} \Bigl) \cdot \Pi_1(s)
\quad (B_1^\prime \geq 0),
\end{equation*}
where
\begin{equation*} \label{18}
\Pi_1(s) = \prod_{0 \not= \beta \in {\Bbb R}} \Bigl( 1 -\frac{s}{\beta} \Bigr)
\prod_{{\rho=\beta+i\gamma}\atop{\beta<1/2,\, \gamma >0}} \left[ \Bigl( 1 -\frac{s}{\rho} \Bigr)\Bigl( 1 -\frac{s}{\overline{\rho}} \Bigr) \right].
\end{equation*}
Note that all zeros of $\Pi_1(s)$ lie in $\sigma_0 < \Re(s) <1/2$ for some $\sigma_0$.
We have
\begin{equation} \label{28}
Z_1(s) = g_1(s) \cdot \Bigl( 1-\frac{g_1(1-s)}{g_1(s)} \Bigr) \quad (g_1(s)=f_1(s)\cdot \chi(2s)).
\end{equation}
and
\begin{equation} \label{29}
\Bigl| \frac{g_1(1-s)}{g_1(s)} \Bigr|
= e^{B_1^\prime (1-2\sigma)}
\cdot \Bigl| \frac{\Pi_1(1-s)}{\Pi_1(s)} \Bigr|
\cdot \Bigl| \frac{s-1+\rho_0}{s-\rho_0}
\cdot  \frac{s-1+\overline{\rho}_0}{s-\overline{\rho}_0} \Bigr|
\cdot \Bigl| \frac{\chi(2s-1)}{\chi(2s)} \Bigr|.
\end{equation}
Because $B_1^\prime \geq 0$, we have
\begin{equation} \label{30}
e^{B_1^\prime(1-2\sigma)} \leq 1 \quad (\Re(s)>1/2).
\end{equation}
For the ratio $\Pi_1(1-s)/\Pi_1(s)$ in \eqref{29}, we have
\begin{equation} \label{31}
\Bigl| \frac{\Pi_1(1-s)}{\Pi_1(s)} \Bigr| =
\prod_{{\rho=\beta+i\gamma}\atop{\beta<1/2,\, \gamma >0}}
\left(
\Bigl| \frac{1-s-\overline{\rho}}{s-\rho} \Bigr| \cdot
\Bigl| \frac{1-s-\rho}{s-\overline{\rho}} \Bigr|
\right) < 1 \quad (\Re(s)>1/2),
\end{equation}
by term-by-term argument as in ~\cite{LS}
by using $\beta<1/2$ and
\[
\Bigl| \frac{1-s-\overline{\rho}}{s-\rho} \Bigr|^2
= 1 - \frac{(2\sigma-1)(1-2\beta)}{(\sigma-\beta)^2 + (t-\gamma)^2},
\]
where $\rho=\beta+i\gamma$ is a zero of $f_1(s)$.
It remains to give an estimate for
\begin{equation}
r_1(s) := \Bigl| \frac{s-1+\rho_0}{s-\rho_0}
\cdot  \frac{s-1+\overline{\rho}_0}{s-\overline{\rho}_0} \Bigr|
\cdot \Bigl| \frac{\chi(2s-1)}{\chi(2s)} \Bigr|.
\end{equation}
To estimate $r_1(s)$, we use the following lemma essentially.
%
%%%%%%%%%%%%%%%%%%%%%%%%%%%%%%%%%%%%%%%%%%%%%%%%%%%%%%%%%%%%%%%%%%%%%%%%%%%%%%%%%%%%
%%
%% lemma 8
%%
%%%%%%%%%%%%%%%%%%%%%%%%%%%%%%%%%%%%%%%%%%%%%%%%%%%%%%%%%%%%%%%%%%%%%%%%%%%%%%%%%%%%
%
\begin{lem}[\cite{La99}] \label{lem_08}
For any real value of $t$  there
exists at least three distinct zeros $\rho= \beta+ i \gamma$ of $\xi(s)$ such that
$0 < \beta \leq 1/2$ and
\begin{equation} \label{35}
|t - \gamma| \leq 22.
\end{equation}
\end{lem}
{\bf Proof.} Suppose $|t| \geq 25$.
Then there exists at least three distinct zeros $\rho=\beta+i\gamma$ of $\xi(s)$
satisfying $0 < \beta \leq 1/2$ and $|t-\gamma| <15.1$
by applying  Lemma 5 in \cite{S} to $t+10.1$ and $t-10.1$
(Lemma 5 in \cite{S} is essentially Lemma 3.5 of \cite{La99}).
For $|t| < 25$, estimate \eqref{35} also holds for three distinct zeros
because $\xi(s)$ has zeros at $s=\pm 14.13, \, \pm 21.02, \, \pm 25.01$.
\hfill $\Box$
\bigskip

\noindent
Using Lemma \ref{lem_08} we show the following.
%
%%%%%%%%%%%%%%%%%%%%%%%%%%%%%%%%%%%%%%%%%%%%%%%%%%%%%%%%%%%%%%%%%%%%%%%%%%%%%%%%%%%%
%%
%% lemma 9
%%
%%%%%%%%%%%%%%%%%%%%%%%%%%%%%%%%%%%%%%%%%%%%%%%%%%%%%%%%%%%%%%%%%%%%%%%%%%%%%%%%%%%%
%
\begin{lem} \label{lem_09}
Let $\rho_0 = \beta_0 + i \gamma_0 \simeq 0.90 + i \cdot 2.09$ be the complex zero of $f_1(s)$ in Lemma \ref{lem_06}.
Let $s = \sigma + it$ with $1/2 < \sigma \leq 20$ and $t \ge 25$.
Then there exists at least two distinct zeros $\rho=\beta+i\gamma$ of $\xi(s)$ such that
$0 < \beta \leq 1/2$, $|t-\gamma| \leq 22$,
\begin{equation} \label{35_1}
\left| \frac{s-1+\overline{\rho_0}}{s-\rho_0} \right| \cdot
\left| \frac{2s-1 - ( 1 - \overline{\rho}) }{2s - \rho} \right|
< 1,
\end{equation}
and
\begin{equation} \label{35_2}
\left| \frac{s-1+\rho_0}{s-\overline{\rho_0}} \right| \cdot
\left| \frac{2s-1 - ( 1 - \overline{\rho}) }{2s - \rho} \right|
< 1.
\end{equation}
%and
%\begin{equation} \label{35_22}
%\left| \frac{s}{s-1} \right| \cdot
%\left| \frac{2s-1 - ( 1 - \overline{\rho}) }{2s - \rho} \right|
%< 1.
%\end{equation}
\end{lem}
{\bf Proof.}
By squaring \eqref{35_1} and \eqref{35_2} we have
\begin{equation} \label{35_3}
\frac{(\sigma+\beta_0-1)^2 + (t \pm \gamma_0)^2}{(\sigma-\beta_0)^2 + (t \pm \gamma_0)^2}
\cdot \frac{(2\sigma+\beta-2)^2 + (t-\gamma)^2}{(2\sigma-\beta)^2 + (t-\gamma)^2} <1.
\end{equation}
To prove Lemma \ref{lem_09} it is sufficient that \eqref{35_3} holds
for $0 < \beta \leq 1/2$, $|t-\gamma| < 22$, $1/2 < \sigma \leq 20$ and $t \geq 25$,
because of Lemma \ref{lem_08}.
To establish \eqref{35_3} in that conditions
it suffices to show that
\begin{equation*}
\frac{(\sigma+\beta_0-1)^2 + (t \pm \gamma_0)^2}{(\sigma-\beta_0)^2 + (t \pm \gamma_0)^2}
\cdot \frac{(2\sigma-\frac{3}{2})^2 + 22^2}{(2\sigma-\frac{1}{2})^2 + 22^2} <1.
\end{equation*}
by a similar reason in the later half of section 4.3 in~\cite{S}.
This inequality is equivalent to
\begin{equation} \label{35_5}
(2\sigma-1) \Bigl( 8(t \pm \gamma_0)^2 - P(\sigma) \Bigr)  > 0,
\end{equation}
where $P(\sigma)=8(4\beta_0-3) \sigma^2 - 8(4\beta_0-3) \sigma - 8 \beta_0^2 + 3890 \beta_0 -  1945$.
Using the value $\beta_0 \simeq 0.90$ we see that $P(\sigma) < 3807$ for $1/2 < \sigma <20$.
On the other hand, using the value $\gamma_0 \simeq 2.09$ we see that
$8(t \pm \gamma_0)^2 > 3872$ for $t \geq 25$
since $|t \pm \gamma_0| = t \pm \gamma_0>22$ for $t \geq 25$.
Hence \eqref{35_5} hold, and it implies \eqref{35_3}.
%
%By a way similar to the above, to establish $\eqref{35_22}$,  it is sufficient to show
%\begin{equation}
%(2\sigma-1) \Bigl( 8t^2 - p(\sigma) \Bigr)  > 0,
%\end{equation}
%where $p(\sigma)=8 \sigma^2 - 8 \sigma + 1937$.
%Because $p(\sigma) < 4977$ for $1/2 < \sigma <20$ and $8*t^2 \geq 5000$
%for $t \geq 25$, we obtain  $\eqref{35_22}$.
\hfill $\Box$
\bigskip

\noindent
Lemma \ref{lem_09} and $\overline{Z_1(s)}=Z_1(\overline{s})$ implies
\begin{equation} \label{36}
|r_1(s)| <1 \quad \text{for} \quad 1/2 < \sigma \leq 20, ~|t| \geq 25
\end{equation}
by taking two distinct zeros of $\xi(s)$ in that region,
since other terms in $r_1(s)$ are estimated as
\begin{equation*}
\Bigl|  \frac{2s-1 -(1-\overline{\rho})}{2s -\rho} \Bigr| < 1 \quad (\Re(s)>1/2),
\end{equation*}
where $\rho$ is a zero of $\xi(s)$.
Estimates \eqref{30}, \eqref{31} and \eqref{36} show that
\begin{equation*}
\Bigl| \frac{g_1(1-s)}{g_1(s)} \Bigr| < 1 \quad \text{for} \quad 1/2 < \sigma \leq 20, ~|t| \geq 25.
\end{equation*}
By \eqref{28} this estimate implies Proposition \ref{prop_03},
because $g_1(s)$ has no zero in the region $1/2 < \sigma \leq 20$, $|t| \geq 25$.
\hfill $\Box$

\subsection{Proof of Proposition \ref{prop_06}}
Let $\rho_0=\beta_0+ i \gamma_0~(\gamma_0>0)$ be the complex zero of $f_2(s)$ in Lemma \ref{lem_07}.
By Proposition \ref{prop_01}, $f_2(s)$ has the factorization
\begin{equation*}
f_2(s) = f_2(0) \, e^{B_2^\prime s} \,
\Bigl(1-\frac{s}{2}\Bigr)
\Bigl( 1-\frac{s}{\rho_0} \Bigl)
\Bigl( 1-\frac{s}{\overline{\rho}_0} \Bigl) \cdot \Pi_2(s)
\quad (B_2^\prime \geq 0),
\end{equation*}
where
\begin{equation*}
\Pi_2(s) = \prod_{0 \not= \beta \in {\Bbb R}} \Bigl( 1 -\frac{s}{\beta} \Bigr)
\prod_{{\rho=\beta+i\gamma}\atop{\beta<1/2,\, \gamma >0}} 
\left[ \Bigl( 1 -\frac{s}{\rho} \Bigr)\Bigl( 1 -\frac{s}{\overline{\rho}} \Bigr) \right].
\end{equation*}
Here all zeros of $\Pi_2(s)$ lie in $\sigma_0 < \Re(s) <1/2$ for some $\sigma_0$.
We have
\begin{equation} \label{28_2}
Z_2(s) = g_2(s) \cdot \Bigl( 1-\frac{g_2(1-s)}{g_2(s)} \Bigr) \quad (g_2(s)=f_2(s)\cdot \chi(2s)).
\end{equation}
and
\begin{equation*}
\Bigl| \frac{g_2(1-s)}{g_2(s)} \Bigr|
= e^{B_2^\prime (1-2\sigma)}
\cdot \Bigl| \frac{\Pi_2(1-s)}{\Pi_2(s)} \Bigr|
\cdot \Bigl| \frac{s}{1-s}  \frac{s-1+\rho_0}{s-\rho_0}
\cdot  \frac{s-1+\overline{\rho}_0}{s-\overline{\rho}_0} \Bigr|
\cdot \Bigl| \frac{\chi(2s-1)}{\chi(2s)} \Bigr|.
\end{equation*}
For $e^{B_2^\prime(1-2\sigma)}$ and $\Pi_2(1-s)/\Pi_2(s)$, we have
\begin{equation} \label{30_2}
e^{B_2^\prime(1-2\sigma)}, \quad \Bigl| \frac{\Pi_2(1-s)}{\Pi_2(s)} \Bigr| < 1 \quad (\Re(s)>1/2)
\end{equation}
by a similar argument as in $f_1(s)$.
It remains to give an estimate for
\begin{equation}
r_2(s) := \Bigl| \frac{s}{1-s}  \frac{s-1+\rho_0}{s-\rho_0}
\cdot  \frac{s-1+\overline{\rho}_0}{s-\overline{\rho}_0} \Bigr|
\cdot \Bigl| \frac{\chi(2s-1)}{\chi(2s)} \Bigr|.
\end{equation}
Using Lemma \ref{lem_08} we show the following:
%
%%%%%%%%%%%%%%%%%%%%%%%%%%%%%%%%%%%%%%%%%%%%%%%%%%%%%%%%%%%%%%%%%%%%%%%%%%%%%%%%%%%%
%%
%% lemma 10
%%
%%%%%%%%%%%%%%%%%%%%%%%%%%%%%%%%%%%%%%%%%%%%%%%%%%%%%%%%%%%%%%%%%%%%%%%%%%%%%%%%%%%%
%
\begin{lem} \label{lem_10}
Let $\rho_0 = \beta_0 + i \gamma_0 \simeq 1.17 + i \cdot 3.43$ be the complex zero of $f_2(s)$ in Lemma \ref{lem_07}.
Let $s = \sigma + it$ with $1/2 < \sigma \leq 20$ and $t \ge 36$.
Then there exists at least three distinct zeros $\rho=\beta+i\gamma$ of $\xi(s)$ such that
$0 < \beta \leq 1/2$, $|t-\gamma| \leq 22$,
\begin{equation} \label{35_1_2}
\left| \frac{s-1+\overline{\rho_0}}{s-\rho_0} \right| \cdot
\left| \frac{2s-1 - ( 1 - \overline{\rho}) }{2s - \rho} \right|
< 1,
\end{equation}
\begin{equation} \label{35_2_2}
\left| \frac{s-1+\rho_0}{s-\overline{\rho_0}} \right| \cdot
\left| \frac{2s-1 - ( 1 - \overline{\rho}) }{2s - \rho} \right|
< 1,
\end{equation}
and
\begin{equation} \label{35_22_2}
\left| \frac{s}{s-1} \right| \cdot
\left| \frac{2s-1 - ( 1 - \overline{\rho}) }{2s - \rho} \right|
< 1.
\end{equation}
\end{lem}
{\bf Proof.}
By squaring \eqref{35_1_2} and \eqref{35_2_2} we have
\begin{equation} \label{35_3_2}
\frac{(\sigma+\beta_0-1)^2 + (t \pm \gamma_0)^2}{(\sigma-\beta_0)^2 + (t \pm \gamma_0)^2}
\cdot \frac{(2\sigma+\beta-2)^2 + (t-\gamma)^2}{(2\sigma-\beta)^2 + (t-\gamma)^2} <1.
\end{equation}
To prove Lemma \ref{lem_10} it is sufficient that \eqref{35_3_2} holds
for $0 < \beta \leq 1/2$, $|t-\gamma| < 22$, $1/2 < \sigma \leq 20$ and $t \geq 25$,
because of Lemma \ref{lem_08}.
To establish \eqref{35_3_2} in that conditions
it suffices to show that
\begin{equation*}
\frac{(\sigma+\beta_0-1)^2 + (t \pm \gamma_0)^2}{(\sigma-\beta_0)^2 + (t \pm \gamma_0)^2}
\cdot \frac{(2\sigma-\frac{3}{2})^2 + 22^2}{(2\sigma-\frac{1}{2})^2 + 22^2} <1.
\end{equation*}
This inequality is equivalent to
\begin{equation} \label{35_5_2}
(2\sigma-1) \Bigl( 8(t \pm \gamma_0)^2 - P(\sigma) \Bigr)  > 0,
\end{equation}
where $P(\sigma)=8(4\beta_0-3) \sigma^2 - 8(4\beta_0-3) \sigma - 8 \beta_0^2 + 3890 \beta_0 -  1945$.
Using the value $\beta_0 \simeq 1.17$ we see that $P(\sigma) < 7777$ for $1/2 < \sigma <20$.
On the other hand, using the value $\gamma_0 \simeq 3.43$ we see that
$8(t \pm \gamma_0)^2 > 8192$ for $t \geq 36$
since $|t \pm \gamma_0| = t \pm \gamma_0 > 32$ for $t \geq 36$.
Hence \eqref{35_5_2} hold, and it implies \eqref{35_3_2}.

By a way similar to the above, to establish $\eqref{35_22_2}$,  it is sufficient to show
\begin{equation*}
(2\sigma-1) \Bigl( 8t^2 - p(\sigma) \Bigr)  > 0,
\end{equation*}
where $p(\sigma)=8 \sigma^2 - 8 \sigma + 1937$.
Because $p(\sigma) < 4977$ for $1/2 < \sigma <20$ and $8t^2 \geq 5000$
for $t \geq 25$, we obtain  $\eqref{35_22_2}$.
\hfill $\Box$
\bigskip

\noindent
Lemma \ref{lem_10} and $\overline{Z_2(s)}=Z_2(\overline{s})$ implies
\begin{equation} \label{36_2}
|r_2(s)| <1 \quad \text{for} \quad 1/2 < \sigma \leq 20, ~|t| \geq 25
\end{equation}
by taking three distinct zeros of $\xi(s)$ in that region.
Estimates \eqref{30_2} and \eqref{36_2} show that
\begin{equation}
\Bigl| \frac{g_2(1-s)}{g_2(s)} \Bigr| < 1
\end{equation}
for $1/2 < \sigma \leq 20$, $|t| \geq 36$.
By \eqref{28_2} this estimate implies Proposition \ref{prop_06},
because $g_2(s)$ has no zero in the region $1/2 < \sigma \leq 20$, $|t| \geq 36$.
\hfill $\Box$

\subsection{Proof of Proposition \ref{prop_04} and \ref{prop_07}}
Because the region $1/2 < \sigma \leq 20$, $|t| \leq 25$ or $36$ is finite,
we can check the assertions of Proposition \ref{prop_04} and Proposition \ref{prop_07} by using the help of computer
as in the proof of Lemma \ref{lem_06} and Lemma \ref{lem_07}. \hfill $\Box$
%
%%%%%%%%%%%%%%%%%%%%%%%%%%%%%%%%%%%%%%%%%%%%%%%%%%%%%%%%%%%%%%%%%%%%%%%%%%%%%%%%%%%%%%%%%%%%%%%%%%%%%%%%%%%%%%%%%%%%%%%%%%
%%
%% section 7
%%
%%%%%%%%%%%%%%%%%%%%%%%%%%%%%%%%%%%%%%%%%%%%%%%%%%%%%%%%%%%%%%%%%%%%%%%%%%%%%%%%%%%%%%%%%%%%%%%%%%%%%%%%%%%%%%%%%%%%%%%%%%
%
\section{Proof of Lemma \ref{lem_01}}
We prove the lemma only if $F(s)$ has genus one,
since if $F(s)$ has genus zero it is easily proved by a way similar to the case of genus one.
The genus one assumption is equivalent to the Hadamard product factorization
\begin{equation}
F(s) = e^{A+Bs} s^m
\prod_{\rho} \Bigl( 1 - \frac{s}{\rho} \,\Bigr) \exp(s/\rho) \quad (m \in {\Bbb Z}_{\geq 0})
\end{equation}
converges absolutely and uniformly
on any compact subsets of $\Bbb C$.
That is also equivalent to $\sum_{\rho} |\rho|^{-2} < \infty.$
Assumption (i) implies the symmetry of the set of zeros under the conjugation $\rho \mapsto \overline{\rho}$.
It follows that the set of zeros $\rho= \beta + i \gamma$,
counted with multiplicity,
is partitioned into  blocks $B(\rho)$ comprising
$\{ \rho, \overline{\rho} \}$ if $\gamma >0$
and
$\{\rho\}$ if $\beta \not=0$ and $\gamma=0$.
Each block is labeled with the unique zero
in it having $\gamma \geq 0$.
Using assumption (ii), we show
\begin{equation} \label{49}
F(s) = s^m e^{A+ B^\prime s}
\prod_{B(\rho)}
\left( \prod_{\rho \in B(\rho)} \Bigl(\, 1 - \frac{s}{\rho} \,\Bigr) \right)
\end{equation}
where the outer product on the right-hand side
converges absolutely and uniformly on any compact
subsets of $\Bbb C$. This assertion holds because the block convergence
factors $\exp(c(B(\rho))s)$ are given by $c(B(\rho)) = 2\beta |\rho|^{-2}$
for $\gamma >0$. Assumption (ii) implies $ |\beta -1/2| < \sigma_0$.
Hence
\[
\sum_{B(\rho)} |c(B(\rho))| \leq
  \sum_{0 \not= \rho :\, {\rm real}} |\rho|^{-1}
+ (2\sigma_0 +1) \sum_{\rho} |\rho|^{-2} < \infty.
\]
Thus the convergence factors $\exp(c(B(\rho))s)$ can be pulled out of the product.
Hence we have $\eqref{49}$ with
\begin{equation}\label{le08}
B^\prime = B + \sum_{B(\rho)}c(B(\rho)).
\end{equation}
Using assumption (iii), (iv) and (v) we show
\begin{equation}\label{le09}
B^\prime  \geq 0.
\end{equation}
By $\eqref{le03}$ in assumption (v) we have
\begin{equation}\label{le13}
{\Bbb R} \ni \log \Bigl( \frac{F(1-\sigma)}{F(\sigma)} \Bigr) \to -\infty
\quad \text{as} \quad \sigma \to +\infty.
\end{equation}
Using $\eqref{49}$ we have
\[
\frac{F(1-\sigma)}{F(\sigma)}
 = e^{B^\prime(1-2\sigma)}
\Bigl( \frac{\sigma-1}{\sigma} \Bigr)^m
\prod_{\rho=\beta \in {\Bbb R}}
\frac{\sigma - 1 + \beta} { \sigma - \beta }
\prod_{{\rho=\beta+i\gamma}\atop{\gamma>0}}
\frac{(\sigma - 1 + \beta)^2 + \gamma^2} { (\sigma - \beta)^2 + \gamma^2 }.
\]
Thus
\begin{equation}\label{le14}
\aligned
\log  \Bigl( \frac{F(1-\sigma)}{F(\sigma)} \Bigr)
& = B^\prime(1-2\sigma)
  + m \log \Bigl( 1-  \frac{1}{\sigma} \Bigr)
  + \sum_{\rho=\beta \in {\Bbb R}}
  \log\Bigl(1 -  \frac{1 - 2\beta} { \sigma - \beta } \Bigr)  \\
& \quad
  +\sum_{{\rho=\beta+i\gamma}\atop{\gamma>0}}
  \log\Bigl( 1- \frac{(1-2\beta)(2\sigma-1)} { (\sigma - \beta)^2 + \gamma^2 }  \,\Bigr).
\endaligned
\end{equation}
Note that
\begin{equation} \label{54}
\log\Bigl( 1- \frac{(1-2\beta)(2\sigma-1)} { (\sigma - \beta)^2 + \gamma^2 } \Bigr) <0
\quad \text{for} \quad \sigma>1/2
\end{equation}
if $\beta<1/2$, and
\[
\log \Bigl( 1-  \frac{1}{\sigma} \Bigr),~
\log\Bigl(1 -  \frac{1 - 2\beta} { \sigma - \beta } \Bigr),~
\log\Bigl( 1- \frac{(1-2\beta)(2\sigma-1)} { (\sigma - \beta)^2 + \gamma^2 } \Bigr)
\to 0 \quad \text{as} \quad \sigma \to +\infty
\]
for any fixed $\rho=\beta+i\gamma$.
By assumption (iii), \eqref{54} holds except for finitely many zeros.
Hence if we suppose $B^\prime <0$, $\eqref{le13}$ and $\eqref{le14}$ implies
\begin{equation}\label{le15}
\left\vert
\sum_{{\rho=\beta+i\gamma}\atop{\gamma>0}}
\log\Bigl( 1- \frac{(1-2\beta)(2\sigma-1)} { (\sigma - \beta)^2 + \gamma^2 }  \,\Bigr)
\right\vert \geq 2 |B^\prime| \sigma
\end{equation}
for large $\sigma>1/2$,
because the number of real zeros is also finite by assumption (ii) and (iii).
On the other hand, for large $\sigma >1/2$, we have
\begin{equation*}
\aligned
\left\vert
\sum_{{\rho=\beta+i\gamma}\atop{\gamma>0}}
\log\Bigl( 1- \frac{(1-2\beta)(2\sigma-1)} { (\sigma - \beta)^2 + \gamma^2 }  \,\Bigr)
\right\vert
& \leq \left\vert
\sum_{{\rho=\beta+i\gamma}\atop{\gamma>0}}
\log\Bigl( 1- \frac{(1-2\sigma_0)(2\sigma-1)} { (\sigma - 1/2)^2 + \gamma^2 }  \,\Bigr)
\right\vert \\
& \ll (2\sigma-1) \sum_{{\rho=\beta+i\gamma}\atop{\gamma>0}}
\frac{1} { (\sigma - 1/2)^2 + \gamma^2 }.
\endaligned
\end{equation*}
The sum in the right-hand side can be written as the Stieltjes integral
\[
\int_{\gamma_0}^{\infty} \frac{dN(t)}{(\sigma-1/2)^2 + t^2}.
\]
Using $\eqref{le02}$ in (iv) we have
\[
\int_{\gamma_0}^{\infty} \frac{dN(t)}{(\sigma-1/2)^2 + t^2}
\ll \int_{\gamma_0}^{\infty} \frac{(\log t) \, dt}{(\sigma-1/2)^2 + t^2}
\ll \frac{\log (\sigma + \gamma_0)}{\sigma-1/2}.
\]
Hence we obtain
\begin{equation}\label{le16}
\left\vert
\sum_{{\rho=\beta+i\gamma}\atop{\gamma>0}}
\log\Bigl( 1- \frac{(1-2\beta)(2\sigma-1)} { (\sigma - \beta)^2 + \gamma^2 }  \,\Bigr)
\right\vert \ll \log(\sigma + \gamma_0)
\end{equation}
for large $\sigma>1/2$.
This contradict $\eqref{le15}$.
Thus $\eqref{le09}$ holds.
\hfill $\Box$
%
%%%%%%%%%%%%%%%%%%%%%%%%%%%%%%%%%%%%%%%%%%%%%%%%%%%%%%%%%%%%%%%%%%%%%%%%%%%%%%%%%%%%%%%%%%%%%%%%%%%%%%%%%%%%%%%%%%%%%%%%%%
%%
%% references
%%
%%%%%%%%%%%%%%%%%%%%%%%%%%%%%%%%%%%%%%%%%%%%%%%%%%%%%%%%%%%%%%%%%%%%%%%%%%%%%%%%%%%%%%%%%%%%%%%%%%%%%%%%%%%%%%%%%%%%%%%%%%
%

\bigskip

\noindent
Masatoshi Suzuki \\
Department of Mathematics \\
Rikkyo University \\
Nishi-Ikebukuro, Toshima-ku \\
Tokyo 171-8501, Japan \\
\texttt{suzuki@rkmath.rikkyo.ac.jp}
\bigskip

\noindent
Lin Weng \\
Graduate School of Mathematics \\
Kyushu University \\
6-10-1, Hakozaki, Higashi-ku \\
Fukuoka 812-8581, Japan \\
\texttt{weng@math.kyushu-u.ac.jp}
\smallskip

\noindent
{\it and}
\smallskip

\noindent
Chennai Mathematical Institute \\
Plot H1, SIPCOT IT Park \\
Padur PO, Siruseri 603103, India

\end{document}